\documentclass[lettersize,journal]{IEEEtran}
\usepackage{amsmath,amsfonts}
\usepackage{algorithmic}
\usepackage{algorithm}
\usepackage{array}
\usepackage[caption=false,font=normalsize,labelfont=sf,textfont=sf]{subfig}
\usepackage{textcomp}
\usepackage{stfloats}
\usepackage{url}
\usepackage{verbatim}
\usepackage{graphicx}
\usepackage{cite}
\usepackage{array}
\usepackage{multirow}
\usepackage{tabularx}
\usepackage{graphicx}
\usepackage{booktabs}
\usepackage{makecell}

\usepackage{color}

\begin{document}

\title{Departure time choice user equilibrium for public transport demand management}


\author{Xia Zhou, Zhenliang Ma,~\IEEEmembership{Member,~IEEE}, Mark Wallace, Daniel D. Harabor 
\thanks{This paper was produced by the IEEE Publication Technology Group. They are in Piscataway, NJ.}
\thanks{Manuscript received February 11, 2025; \textit{(Corresponding author: Zhenliang Ma.)}}}

\markboth{Journal of \LaTeX\ Class Files,~Vol.~14, No.~8, February~2025}%
{Shell \MakeLowercase{\textit{et al.}}: A Sample Article Using IEEEtran.cls for IEEE Journals}


\maketitle

\begin{abstract}
Departure time management is an efficient way in addressing the peak-hour crowding in public transport by reducing the temporal imbalance between service supply and travel demand. From the demand management perspective, the problem is to determine an equilibrium distribution of departure times for which no user can reduce their generalized cost by changing their departure times unilaterally. This study introduces the departure time choice user equilibrium problem in public transport (DTUE-PT) for multi-line, schedule-based networks with hard train capacity constraints. We model the DTUE-PT problem as a Non-linear Mathematical Program problem (NMP) (minimizing the system gap) with a simulation model describing the complex system dynamics and passenger interactions. We develop an efficient, adaptive gap-based descent direction (AdaGDD) solution algorithm to solve the NMP problem. We validate the methodology on a multi-line public transport network with transfers by comparing with classical public transport assignment benchmark models, including Method of Successive Average (MSA) and day-to-day learning methods. The results show that the model can achieve a system gap ratio (the solution gap relative to the ideal least cost of an origin-destination option) of 0.1926, which significantly improves the solution performance from day-to-day learning (85\%) and MSA (76\%) algorithms. The sensitivity analysis highlights the solution stability of AdaGDD method over initial solution settings. The potential use of DTUE-PT model is demonstrated for evaluating the network design of Hong Kong mass transit railway network and can be easily extended to incorporate the route choice.

\end{abstract}


\section{Introduction}
\IEEEPARstart{M}{anaging} peak-hour congestion has become a critical challenge for transportation policymakers and researchers. Peak-hour congestion arises from an imbalance between the limited supply of services and the high travel demand in certain places at certain times. The modeling and analysis of peak-hour congestion have a rich history, starting with the pioneering work from \cite{ref1}. Vickrey \cite{ref1} proposed the bottleneck model in which each user chooses their departure time to minimize their own generalized travel cost. The problem is to determine a dynamic equilibrium distribution of departure times, known as the Departure Time Choice User Equilibrium (DTUE). At the DTUE, no user can decrease their perceived travel costs by unilaterally changing their departure times (known as Wardrop's first principle \cite{ref25}). Insights into the characteristics of the DTUE problem contribute significantly to the evaluation and formulation of effective transport policies to alleviate peak-hour congestion \cite{ref2}. For example, if incentives were proposed, the DTUE model could be used to predict the potential impact of these incentives on users' behavior.

The DTUE problem has been widely studied in car traffic \cite{ref3}, but no such problem is studied in the public transport domain (named DTUE-PT in the paper). Unlike the open system of the car traffic, the closed system of public transport is constrained by predefined routes and schedules \cite{ref4}. These bring additional modeling challenges including, for example, transfer times between service lines and waiting times and denied boarding times at the boarding station due to constrained vehicle capacity. Importantly, the discrete, nonlinear nature of system dynamics and travel costs makes the DTUE-PT problem more complex compared to the continuous-flow car traffic system. Therefore, DTUE models in car traffic are not easily applied or extended to public transport.

Although no exact DTUE-PT problem was studied in the public transport literature, Table \ref{tab:table-1} lists the relevant UE models reported and details their problem focus, mathematical objectives, solution algorithms, and modeling details (hard capacity and multi-line network). Mathematically, various objective functions of the UE model are proposed, including System Gap Minimization \cite{ref6, ref13}, Stochastic User Equilibrium Equation \cite{ref8, ref18, ref20}, Variational Inequality \cite{ref12, ref16, ref19, ref22}, Integral UE Formulation \cite{ref9, ref11}, and Deterministic User Equilibrium Equations \cite{ref15, ref21}. 

However, it remains unclear whether these mathematical models satisfy Wardrop's first principle (UE assignment) \cite{ref25}. Also, no study examined its solution with respect to the UE assignment principle which is important for the UE model benchmark and ultimately for decision making practices.

\begin{table*}[ht]
    \centering
    \caption{Literature review: UE problems in public transport}
    \begin{tabular}{p{2.4cm} p{2.4cm} p{2.0cm} p{2.0cm} p{2.0cm} p{2.8cm}} 
        \toprule
        \textbf{Publications} & \textbf{Focus} & \textbf{Objective} & \textbf{H-capacity} & \textbf{M-line} &
        \textbf{Algorithm}\\
        \midrule
        \cite{ref5} & Route & - & \texttimes & \checkmark & DTD Learning\\
        \cite{ref6} & Route & SGM & \checkmark & \texttimes & MSA\\
        \cite{ref7} & Route \& DT & - & \texttimes & \texttimes & MSA\\
        \cite{ref8} & Mode \& DT & SUE & \checkmark & \texttimes & MSA\\
        \cite{ref9} & Route \& DT & I-UE & \checkmark & \texttimes & MSA\\
        \cite{ref10} & DT & Theoretical & \checkmark & \texttimes & -\\
        \cite{ref11} & Route & I-UE & \texttimes & \checkmark & MSA\\
        \cite{ref12} & Strategy & VI & \checkmark & \texttimes & MSA\\
        \cite{ref13} & Mode & SGM &  \checkmark (partially) & \checkmark & MSA\\
        \cite{ref14} & Station & Theoretical & \texttimes & \texttimes & -\\           
        \cite{ref15} & DT & DUE & \checkmark & \texttimes & Theoretical\\
        \cite{ref16} & Route & VI & \texttimes & \checkmark & Newton-type\\             
        \cite{ref17} & Fare scheme  & Theoretical & \texttimes & \texttimes & -\\             
        \cite{ref18} & Route & SUE & \texttimes & \texttimes & MSA\\        
        \cite{ref19} & DT & VI & \texttimes & \texttimes & Route-swapping\\
        \cite{ref20} & Route & SUE & \texttimes & \texttimes & Branch \& bound\\        
        \cite{ref21} & DT & DUE & \texttimes & \texttimes & Theoretical\\
        \cite{ref22} & Route & VI & \checkmark & \texttimes & Diagonalized\\    
        \cite{ref36} & Route & - & \texttimes & \checkmark & DTD Learning\\
        \cite{ref37} & Mode\& Route & - & \checkmark & \texttimes & DTD Learning\\
        \cite{ref38} & DT\& Route & - & \checkmark & \texttimes & DTD Learning\\
        \cite{ref40} & Route & Fixed-point & \checkmark & \texttimes & MSA\\
        \textbf{This study} &  DT\& Route & SGM & \checkmark & \checkmark & AdaGDD \\ 
        \multicolumn{6}{l}{} \\
        \multicolumn{6}{l}
        {\textbf{H-capacity:} Hard capacity;  \textbf{M-line: } Multi-line networks; \textbf{DT:} Departure time (choice)}  \\
        \multicolumn{6}{l}
        {\textbf{VI: } Variational inequality; \textbf{SUE: } Stochastic user equilibrium equation;}  \\
        \multicolumn{6}{l}{  \textbf{DUE: }Deterministic user equilibrium equation; \textbf{I-UE: } Integral UE formulation;} \\
        \multicolumn{6}{l}{\textbf{I-UE: } Integral UE Formulation; \textbf{SGM: } System gap minimization; } \\
        \bottomrule
    \end{tabular}
    \label{tab:table-1}
\end{table*}

Model tractability is a major challenge for the DTUE-PT problem, as passengers' travel cost decisions are interconnected, the travel cost function is nonlinear, and the system is constrained by `hard' capacity (leading to denied boarding passengers) \cite{ref2}. Given these, modeling waiting times in public transport remains a challenge when considering hard capacity constraints. Table \ref{tab:table-1}  shows that most studies neglect or use soft capacity constraints, assuming that no passengers are left on the platform when a train departs (squeeze passengers in regardless of the remaining capacity). Furthermore, given the complexity of the problem, many UE studies consider only the single-line network (with no transfer) in public transport. No DTUE-PT study considered both the hard vehicle capacity and the multi-line network with transfers in public transport. Such problem is challenging for its nonlinear, non-tractable modeling of passengers' travel costs given complex passenger interactions, transfers, and denied boarding at stations.

Given the intractable problem, heuristic algorithms are commonly used to solve the optimization model for the UE assignment. Typical algorithms of UE problems include (see Table \ref{tab:table-1}): (a) the Method of Successive Averages (MSA) \cite{ref6}-\cite{ref9}, which iteratively updates flow assignments toward equilibrium by averaging solutions across iterations; (b) the Day-to-Day (DTD) Learning method (or day-to-day dynamical approach, day-to-day evolution) \cite{ref36}-\cite{ref38}, which simulates travelers' learning behavior over successive days, gradually converging toward equilibrium by mimicking users' adaptation to service conditions. 
The main disadvantage of the MSA is that it treats all the path (or option choice) equally \cite{ref13}. For a given origin-destination (OD) pair, the relative shift in flows to the current `best' path from a slightly `bad' path is the same as from a severely `bad' path. The DTD learning method relies on predetermined flow updating strategies based on passengers' previous actual travel costs and predicted travel costs. However, it neglects the gap between `good' and `bad' paths. Lu et al. \cite{ref31} proposed a gap-based path-swapping descent direction method for the UE traffic assignment problem, where the probability of a vehicle switching paths is proportional to the relative gap between the experienced path and the best path. This method avoids the need to compute the gradient of the objective function when determining descent directions, making it particularly valuable for large-scale user equilibrium problems. To the best of our knowledge, no study has applied this method to public transport in the context of the DTUE problem.

To address the identified research gaps, the paper introduces the DTUE-PT problem with Wardrop's first principle (UE assignment), considering simultaneously hard train capacity constraints and multi-line networks. We formulate the DTUE-PT problem as a nonlinear mathematical program (NMP) with a simulation model describing the system dynamics and interactions with passengers. The objective function is minimizing the system gap of travel cost between non-optimal options and the optimal option (with the least average cost), and the decision variables are the distributions of the passenger flow over different departure times. We develop an adaptive gap-based descent direction algorithm (AdaGDD) to find the local optimal solution to the NMP problem. We validate the methodology using a multi-line public transport network with transfers, by comparing it to public transport assignment benchmark methods: MSA and DTD learning. The main contributions are: 
\begin{itemize}
\item Introduce the DTUE-PT problem with Wardrop's first principle (UE assignment), considering simultaneously hard train capacity constraints and multi-line networks.
\item Formulate the DTUE-PT problem as a nonlinear mathematical program (NMP) with a simulation model describing the system dynamics and interactions with passengers, and an objective function minimizing the system gap of average travel cost between non-optimal options and the optimal option (with the least average cost).
\item Develop an AdaGDD algorithm that utilizes an optimized step-size determination and updating mechanism to solve the NMP.
\item Validate the methodology using a multi-line network by comparing it with the MSA and DTD learning methods and performing sensitivity analysis. 
\item The proposed method is demonstrated by evaluating the network design of the Hong Kong Mass Transit Railway (MTR) and its extensions, incorporating route choice.
\end{itemize}

The remainder of the paper is organized as follows. Section \ref{sec2}  formulates the DTUE-PT problem, and its solution algorithm is developed in Section \ref{sec:solution-algorithm}. Section \ref{sec:case-study} validates the proposed methodology using a multi-line public transport network, conducts the sensitivity analysis of the AdaGDD algorithm, and demonstrates the use of the DTUE-PT model. The last section summarizes the main findings and future research.

\section{Problem Formulation}\label{sec2}
We consider a public transport network $G = (N, L)$, where $N$ is the set of stations and $L$ is the set of links. The network is composed of distinct elements: a set of origin stations ($O \subset N$), a set of destination stations ($D \subset N$), and a set of transfer stations ($W \subset N$).

The set $K$ represents OD pairs involved in the public transport system. Each origin station ($o \in O$) offers train set $T_o$. The demand between OD pair $k \in K$ is represented by $Q_k$, which upon assignment to the public transport system, generates departure time flows $q_k(t) \in \mathbb{Z^+}$, where $t \in T_{o(k)}$. $T_{o(k)}$ is the departure time set of trains at origin station $o(k)$ (the origin station of OD pair $k$). 

The DTUE-PT problem is: given a multi-line, schedule-based public transport network \cite{ref39} with hard train capacity constraints (see Eq. \ref{eq:88}), determine the equilibrium distribution of departure times for which no user can reduce their generalized costs by changing their departure times unilaterally. The following sub-sections formulate the DTUE-PT model. In the model, we assume that: 1) All users have the same desired arrival time at destinations (i.e., work time); 2) The train schedules are fixed with no randomness (i.e., on-time), and 3) Waiting passengers on the platform are mingled perfectly and board trains based on the first come first serve (FCFS) principle. Key notations used in this study are listed in Table \ref{tab2}.

\subsection{Generalized travel costs in PT}
In this section, each OD pair is assigned a single route ($r$), simplifying the analysis. The model is later extended to incorporate simultaneous departure time and route choices in Section \ref{sec:result4}.  Here, users only choose their departure time ($t$) at the origin. Each departure time on that OD corresponds to the time a train departs from the origin (the user's departure time is equivalent to the same departure time of a train at the same origin).  We focus on the departure time choice problem here and assume the fixed or predetermined known path set per OD pair. Therefore, it is a discrete time choice problem. This subsection summarizes the evaluation of individual travel costs for all users $\mathbf{\textit{C}} = \{c_k^i(t), k \in K, t \in T_{o(k)}, i \in I_k(t)\}$ ($I_k(t)$ is user group who choose departure time $t$ within OD pair $k$). The set of decision variable is the demand flow distribution $\mathbf{\textit{Q}} = \{q_k(t), k \in K, t \in T_{o(k)}\}$. 
\begin{table}[H]
    \begin{minipage}[b]{0.45\textwidth}
        \centering
        \caption{Notations - System Parameters}
        \label{tab2}
        \begin{tabularx}{\columnwidth}{cX X}
            \toprule          
            \multicolumn{2}{l}{\textbf{Network Parameters}}\\ 
            \midrule
            $N$ & Set of stations\\
            $L$ & Set of links \\
            $O$ & Set of origins\\
            $D$ & Set of destinations\\
            $W$ & Set of transfer stations\\
            $K$ & Set of OD pairs ($K\subseteq O \times D$)\\
            $T_o$ & Set of train runs at origin station $o$ (default: 100 trains)\\ 
            $\Omega_o$ & Capacity of  trains departing from $o$ (default: 230 users)\\
            $\pi_o$ & Frequency of trains departing from $o$ (default: 5 mins)\\
            \midrule
            \multicolumn{2}{l}{\textbf{Model Parameters}}\\ 
            \midrule
            $t^*$ & Expected work start time (default: 9:00 AM)\\
            $\alpha$ & Unit cost values of waiting time (default: 10)\\
            $\beta$ & Unit cost values of early arrival penalty (default: 1)\\
            $\gamma$ & Unit cost values of late arrival penalty (default: 10)\\            
            \midrule
            \multicolumn{2}{l}{\textbf{Model Variables}}\\ 
            \midrule
            $q_k(t)$ &  Departure time flow distribution\\
            $c_k^i (t)$ & Individual travel cost of the $ith$ user in $(k, t)$, who plans to take $t$, $i \in I_k(t)$; ($I_k(t)$: the user group ($k, t$))\\
            $w_k^i(t)$ & Waiting time of the user ($k, t, i$)\\
            $\varphi(a_{k}^i (t))$ & Early or late penalty cost of the user ($k, t, i$)\\
            $a_k^i(t)$ & Arrival time of the user ($k, t, i$)\\
            $v_{k}^i (t)$ & In-vehicle time of the user ($k, t, i$)\\
            $g_o(t)$ & The number of waiting users at $o$ for train $t$, ($g_b(t)$: at transfer station $b$, $b \in W$)\\
            $h_o(t)$ & The number of denied boarding users by train $t$ at $o$, ($h_o^k(t)$: for OD pair $k$; $h_b(t)$: at transfer station $b$)\\
            $f_o(t)$ & The number of successfully boarded users on train $t$ at $o$, ($f_o^k(t)$ for OD pair $k$, $f_b(t)$ at transfer station $b$)\\
            $z_o(t)$ & The number of in-vehicle users when train $t$ leaves $o$, ($z_b(t)$: at transfer station $b$)\\  
            $m_b(t)$ & The number of alighting users when train $t$ arrives at $b$ ($m_b^k(t)$ for OD pair $k$; $m_d(t)$ at destination $d$, $d \in D$)\\ 
            $x_b(t)$ & The available capacity when train $t$ completes dropping off users at $b$\\
            $y_b(t)$ & The number of current users in the train at $b$\\
            $\mathbf{\textit{Q}}$ & Demand flow distribution $\{q_k(t), k \in K, t \in T_{o(k)}\}$\\
            $\mathbf{\textit{M}}$ & Arrival flow distribution $\{m_d^k(t), d \in D,k \in O \times \{d\}, t \in T_{o(k)}\}$\\
            $\mathbf{\textit{A}}$ & Arrival time distribution $\{a_k^i(t), k \in O \times \{d\}, t \in T_{o(k)}, i \in I_k(t)\}$\\
            $\mathbf{\textit{C}}$ & Individual travel cost distributions $\{c_k^i(t), k \in O \times \{d\}, t \in T_{o(k)}, i \in I_k(t)\}$ ($C_k(t)$: the average cost for option ($k, t$); $C_k^0(t)$: the free flow cost)\\
            $C_k^*$ & Optimal option cost (least average cost) within $k$\\
            $\mathbf{\sigma}$ & Actual step size distribution
            $\{\sigma_k(t), k \in K, t \in T_{o}\}$\\
            $\mathbf{\partial}$ & OD weighted ratio distribution
            $\{\partial_k, k \in K\}$\\
            $\mathbf{\varphi}$ & Weighted time option distribution
            $\{\varphi_k(t), k \in K, t \in T_{o}\}$\\
            $\theta$ & System step size\\
            $\mathbf{\tau}_k$ & Set of non-optimal time options\\
            $\mathbf{\zeta}(\mathbf{\textit{Q}})$ & The value of the objective function with variable $\mathbf{\textit{Q}}$\\
            \bottomrule
        \end{tabularx}
    \end{minipage}

\end{table}
\begin{figure*}[ht]
    \centering
    \includegraphics[width=1\textwidth]{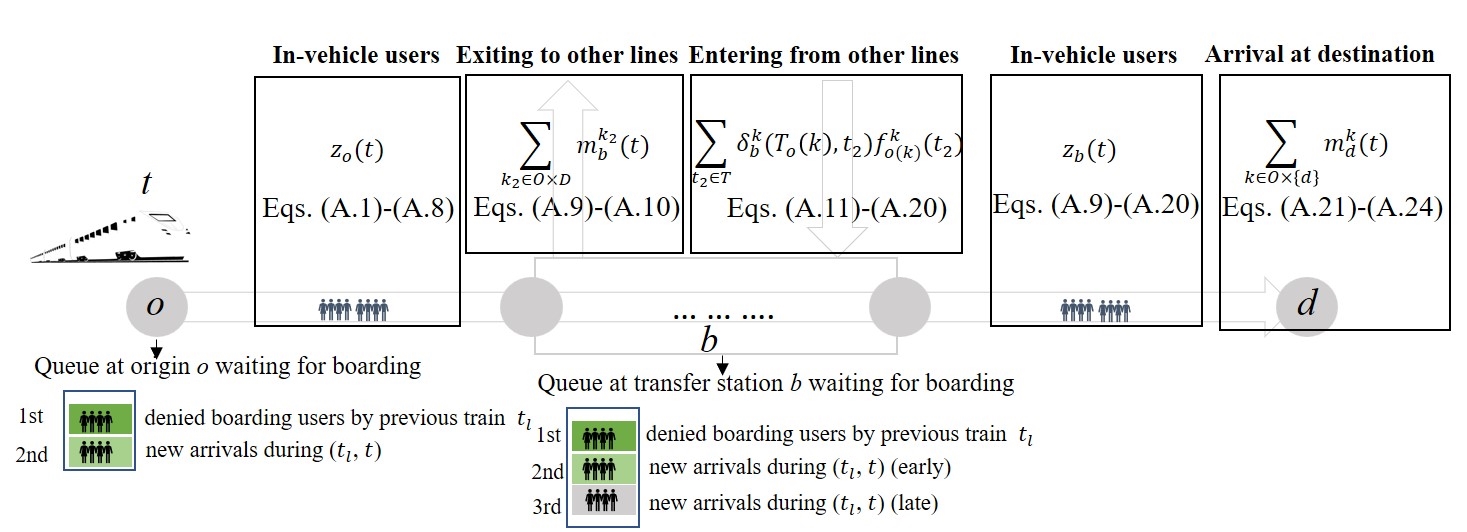}
    \caption{Flow conservation in a congested PT network for evaluating arrival times (see Appendix-A)}
    \label{fig:fig28}
\end{figure*}
Let $Q_k$ be the total demand of OD pair $k$, then
\begin{equation}
    \label{eq:1}
    Q_k = \sum_{t \in T_{o(k)}} q_k(t)
\end{equation}

Eq. (\ref{eq:1}) ensures that the summation of time option flows $q_k(t)$ corresponds to the total OD demand $Q_k$.

The generalized travel cost linked to a time option comprises two elements: the waiting cost and the cost related to schedule delay, whether early or late. The model can also be extended to include in-vehicle cost, fare prices and/or incentives by updating Eq. (\ref{eq:4}), which is valuable for predicting the potential impact of incentives on system performance. The perceived individual travel cost is denoted as $c_k^i (t)$:
    \begin{equation}
    c_k^i (t) = \alpha w_{k}^i(t) + \varphi(a_{k}^i (t))
    \label{eq:4}
    \end{equation}
where $\alpha$ is the unit waiting time weighting and $w_{k}^i(t)$ is the waiting time. $\varphi(a_{k}^i (t))$ represents the early or late arrival penalty incurred by the user. $w_{k}^i(t)$ is given by:
\begin{equation}
w_{k}^i (t) = a_{k}^i (t) - v_{k}^i (t) - t
\label{eq:5}
\end{equation}
where $a_{k}^i (t)$ is the arrival time, and $v_{k}^i (t)$ is the in-vehicle time. In this study, the in-vehicle times are constant for all time options. $a_{k}^i (t)$ is a complicated variable, which will be discussed in the next section.

The early or late penalty cost in Eq. (\ref{eq:4}) is given by:
\begin{equation}
\varphi\left(a_{k}^i (t)\right) = \begin{cases}
    \beta \left(t_k^* - a_{k}^i (t)\right), &  a_{k}^i (t) < t_k^* \\
    0, & a_{k}^i (t) = t_k^*\\
    \gamma \left(a_{k}^i (t) - t_k^* \right), & a_{k}^i (t) > t_k^*
\end{cases}
\label{eq:6}
\end{equation}
In this context, $\beta$ and $\gamma$ signify weightings corresponding to the early and late arrival times, respectively. Additionally, $t_k^{*}$ represents the user’s preferred work start time, which is set uniformly to 32400 seconds (equivalent to 9:00 AM) for all users throughout this study.

\subsection{Transit assignment and simulation}
To evaluate $a_{k}^i (t)$, we propose a transit assignment model to capture flow conservation in a congested public transport (PT) network, as shown in Appendix-A. Fig. \ref{fig:fig28} displays the flowchart of the assignment model. It models trains ($t \in T$) picking and dropping users at different stations (origin station $o \in O$, transfer station $b \in B$, destination station $d \in D$), which can output arrival flows for different OD pairs $\mathbf{\textit{M}}$. Then we can evaluate arrival times for all users $\mathbf{\textit{A}}$ by mapping $\mathbf{\textit{M}}$ to $\mathbf{\textit{Q}}$ according to the FCFS rule.

There are several discontinuities in boarding/denied boarding and transfers/missed connections in the transit assignment model. In order to capture these discontinuities, we propose a simulation model based on the transit assignment model. The simulation model (network loading) is based on the event-based network performance model structure in \cite{ref29} and \cite{ref30}. It works as follows: the events in the system include the arrival and departure of train runs at stations.  Events are sequenced in time.  At each train departure, the users on the train comprise: (1) the users already on the train before this station, and (2) a fraction of the users waiting at the station platform for this train: the fraction is 1 if there is enough room on the train, otherwise it is the fraction that fills the train to capacity.  At each train arrival, the users leaving the train are: (1) those who have arrived at their destination, and (2) those who are changing to the next train on their routes.  These users proceed to the relevant platform and wait for the next train. The input to the simulation is the departure time distribution at origin stations for all users $\mathbf(\textit{Q})$.  Its output is the arrival time distribution at destination stations for all users ($\mathbf{\textit{A}}$). 

\subsection{User equilibrium of departure time choices}
A DTUE-PT solution adheres to Wardrop’s user equilibrium principle when the perceived travel cost for users across all available time options for an OD pair is equal and minimized. Conversely, the perceived travel cost for unused time options is higher than or equal to that of the utilized options. This principle can be mathematically expressed as follows:
\begin{subequations}\label{eq:30}
    \begin{align}
        &C_k(t)= C_k^* \quad \text{if } q_k(t) > 0, &\forall k \in K, t \in T_{o(k)}
        \label{eq:30a}\\
        &C_k(t) \geq C_k^*, & \forall k \in K, t \in T_{o(k)}
        \label{eq:30b}\\
        &q_k(t) \geq 0, & \forall k \in K, t \in T_{o(k)}
        \label{eq:30c}\\
        &Q_k - \sum_{t \in T_{o(k)}} q_k(t) = 0, &  \forall k \in K\label{eq:30d}
    \end{align}
\end{subequations}
where $C_k(t)$ denotes the average travel cost of departure time option $t$ within OD pair $k$. $C_k^*$ is the minimal average travel cost (optimal perceived travel cost) among $T_{o(k)}$. Constraint (\ref{eq:30a}) ensures that the perceived travel costs of all used time options are the same and equal to the optimal perceived travel cost $C_k^*$. Constraint (\ref{eq:30b}) ensures that the perceived travel costs of all time options are not less than the optimal perceived travel cost. Constraint (\ref{eq:30d}) ensures that the sum of time option flows ($\sum_{t \in T_{o(k)}} q_k(t)$) within OD pair $k$ equals to the total demand of that OD $Q_k$. $C_k(t)$ is given by:
\begin{equation} 
C_k(t) =  \begin{cases} 
    \frac{\displaystyle \sum_{i \in I_k(t)} c_k^i(t)}{\displaystyle q_k(t)} & \text{\textit{if}} q_k(t) > 0\\
    \\
    C_k^{0}(t) & \text{\textit{if}}q_k(t) = 0
\end{cases}
\label{eq:31}
\end{equation}
where $C_k^{0}(t)$ is the free flow cost of the option ($k, t$). $C_k^*$ is given by:
\begin{equation}
C_k^* = \min_{t \in T_{o(k)}} C_k(t), \quad \forall k \in K
\label{eq:32}
\end{equation}

The UE conditions expressed in Eqs. \eqref{eq:30a}-\eqref{eq:30d} can alternatively be reformulated as a non-linear complementarity problem (NCP), as discussed in the literature \cite{ref24}.

\begin{subequations}\label{eq:33}
    \begin{align}
        &\sum_{k \in K} \sum_{t \in T_{o(k)}}\left(C_k(t) - \rho_k\right)q_k(t) = 0
        \label{eq:33a}\\
        &C_k(t) - \rho_k \geq 0, \quad \forall k \in K, t \in T_{o(k)} \label{eq:33b}\\
        &q_k(t) \geq 0, \quad \forall k \in K, t \in T_{o(k)} \label{eq:33c}\\
        &Q_k - \sum_{k \in K} q_k(t) = 0, \quad \forall k \in K \label{eq:33d}
    \end{align}
\end{subequations}

Constraint (\ref{eq:33a}) ensures that the system gap (the sum of the gaps among time options and the `optimal' time option) is zero. Constraints (\ref{eq:33b}) and (\ref{eq:33c}) are same as constraints (\ref{eq:30b}) and (\ref{eq:30c}). Constraint (\ref{eq:33d}) ensures that the sum of time option flows equals the total OD demand, thereby maintaining flow conservation. The aforementioned NCP represents an idealized UE condition where the perceived generalized travel costs among all utilized time options for the same OD pair are identical, resulting in a zero-gap solution. However, practical implementation often deviates from this ideal scenario (NCP) due to interacting OD flows (Proposition 1) within the public transport network. For example, reducing the gap within OD pair O1-D1 may increase the gap within OD pair O2-D1, as illustrated in Fig. \ref{fig:fig4}).

\textbf{Proposition 1.} In the DTUE-PT solution, two public transport routes (corresponding to two OD pairs) that are separated by time on the same train line may conflict with each other if there is a sufficiently large flow on the earlier route. (e.g., shifting users among OD1 to reduce the option cost gap for OD1 may inevitably increase the option cost gap for OD2). The proof example of \textbf{Proposition 1.} can be found in Appendix-C.

From \textbf{Proposition 1}, in complex public transport networks, such as multi-line systems, the public transport network equilibrium problem (NCP) may not have a clear solution. Such UE conditions were formulated as a gap-minimizing mathematical model in road traffic \cite{ref31}, where the probability of a vehicle switching paths is proportional to the relative gap between the best path and the current path. Recently, \cite{ref13} adopted the proposed gap-minimizing method based on \cite{ref31} to the UE problem in public transport. When transferring just one user from a worse departure time to any better departure time on any OD pair fails to improve the gap, the corresponding solution can be considered as the UE solution. This paper formulates a gap-based UE model for the DTUE-PT assignment problem considering simultaneous hard train capacity and multi-line networks. Compared to \cite{ref13}, this paper differs in three key aspects. First, our model can be extended to incorporate simultaneous departure time and route choice, rather than focusing solely on the mode choice problem. Second, we account for hard train capacity constraints instead of soft capacity. Finally, we develop a two-loop AdaGDD method instead of an MSA-based approach. 

We replaced Eq. (\ref{eq:33a}) with an objective function, Eq. (\ref{eq:34a}), aimed at minimizing the system gap, thereby bringing it closer to zero as specified in Eq. (\ref{eq:33a}), as follows:
\begin{subequations}\label{eq:34}
 \begin{align}
        &\min_{\mathbf{\textit{Q}}, \mathbf{\rho}}\zeta(\mathbf{\textit{Q}, \mathbf{\rho}}) = \sum_{k \in K} \sum_{t \in T_{o(k)}}\left(C_k(t) - \rho_k\right)q_k(t) 
        \label{eq:34a}\\
        &C_k(t) - \rho_k \geq 0, \quad \forall k \in K, t \in T_{o(k)}
        \label{eq:34b}\\
        &q_k(t) \geq 0, \quad \forall k \in K, t \in T_{o(k)} \label{eq:34c}\\
        &Q_k -  q_k(t) = 0, \quad \forall k \in K \label{eq:34d}\\
        &\rho_k \geq 0, \quad \forall k \in K
        \label{eq:34e}
    \end{align}
\end{subequations}

The objective function Eq. (\ref{eq:34a}) consists of two control variables: $\mathbf{\textit{Q}}$, and $\mathbf{\rho}$, $\mathbf{\rho}=\{\rho_k,k \in K \}$, aiming to minimize the gap to find the optimal solution. When the optimal solution ($\mathbf{\textit{Q}}^*$, $\mathbf{\rho}*$) is found,  $\rho_k^*$ will be the optimal travel cost of the OD $k$ ($\rho_k^*= \min_{t \in T_{o(k)}}C_k(t),  k \in K$) (refer to \textbf{Proposition 2}). 

\textbf{Proposition 2.} At the optimal solution of the Non-linear Mathematical Program (NMP), the optimal $\rho_k^*$ represents the minimum perceived travel costs, given by ($\rho_k^* = \min_{t \in T_{o(k)}}C_k(t), k \in K$).

\textbf{Proof.} At the optimal solution of the Non-linear Mathematical Program (NMP), both $C_k^*$, and $\rho_k^*$ represent the optimal values. Assuming $\rho_k^*$ is not the minimum value among \{$C_k(t), k \in K, t \in T_{o(k)}$\} based on Eq. \eqref{eq:34a}, it follows that $\rho_k^* < C_k(t)$. This discrepancy implies a positive gap $\eta = \min_{t \in T_{o(k)}}C_k(t) - C_k^*$ for any $k$, $k \in K$. However, replacing $\rho_k^*$ with $\rho_k^* + \eta$ in the objective function Eq. \eqref{eq:34a} leads to a further reduction in its value, contradicting the optimality of $\rho_k^*$. Thus, $\rho_k^* = \min_{t \in T_{o(k)}}C_k(t), k \in K$. 

Therefore, the NMP problem can be treated as an equivalent minimization problem:
\begin{subequations}\label{eq:35}
 \begin{align}
        &\min_{\mathbf{\textit{Q}}}\zeta(\mathbf{\textit{Q}}) = \sum_{k \in K} \sum_{t \in T_{o(k)}} \left(C_k(t) - C_k^*\right)q_k(t) 
        \label{eq:35a}\\
        &q_k(t) \geq 0, \quad \forall k \in K, t \in T_{o(k)} \label{eq:35b}\\
        &Q_k -  \sum_{t \in T_{o(k)}} q_k(t) = 0, \quad \forall k \in K
        \label{eq:35c}
    \end{align}
\end{subequations}

The objective function in Eq. (\ref{eq:35a}) only involves control variable $\mathbf{\textit{Q}}$, and aims to minimize the system gap to find the local optimal solution. Assume the NMP has a current solution $\mathbf{\textit{Q}}^m$, and the next potential solution is $\mathbf{\textit{Q}}^{m+1}$. We can calculate the system gap $\zeta(\mathbf{\textit{Q}}^m)$ using Eqs. (\ref{eq:31}), (\ref{eq:32}) and (\ref{eq:35a}). Then, we search for the descent direction and step size to find the next potential solution $\mathbf{\textit{Q}}^{m+1}$ for $\zeta(\mathbf{\textit{Q}}^{m+1})$ using the same method. The detailed iterative solution algorithm will be explained in Section \ref{sec:solution-algorithm}.

\section{Solution Algorithm}\label{sec:solution-algorithm}
We develop an AdaGDD solution algorithm to solve the DTUE-PT model, as outlined in Algorithm \ref{alg:solutionalgorithm}. The probability of a passenger switching options (departure times) is proportional to the relative gap between the experienced option cost and the cost of the best option. The system gap is guaranteed to decrease after each successful iteration until it can no longer be reduced, even by shifting a single passenger from a  `bad' option to a `good' one for any OD pair.

The AdaGDD algorithm comprises two loops. The outer loop, which implements a system-based flow shifting strategy, aims to find the optimal coarse assignment of time option flows that minimizes the system gap. Based on the `coarse` solution, the inner loop (OD-based flow shifting strategy) is used to refine the time option flow assignment and further minimize the system gap. The inner loop ends when transferring just one user from a non-optimal departure time to any optimal departure time, on any OD, fails to improve the system gap.  Both loops employ the gap-based descent direction method based on the framework of \cite{ref31}, incorporating golden section search method \cite{ref32}. 

\begin{algorithm}
\caption{AdaGDD Solution Algorithm} 
\label{alg:solutionalgorithm}
\begin{algorithmic}[1]
\STATE \textbf{Procedure} {(1) System-based flow shifting loop: }({\textsc{Input}}: $ \mathbf{\textit{Q}} $; {\textsc{Output}}: $ \mathbf{\textit{Q}_\textit{1}^\textit{*}}$ )\\
\FOR{$j \leftarrow 0$ to $N$}
\STATE {\textsc{Evaluation}}: $\mathbf{\textit{M}_\textit{j}} \gets \mathbf{\textit{Q}_\textit{j}} $, via Simulation; $\mathbf{\textit{A}_\textit{j}} \gets \mathbf{\textit{M}_\textit{j}}$; (subscripted $\mathbf{\textit{Q}_\textit{j}}$: $\mathbf{\textit{Q}}$ in the $j$th iteration)
\STATE \hspace{0.5cm} $\mathbf{\textit{C}_\textit{j}} \gets \mathbf{\textit{A}_\textit{j}}$, via Eqs. (\ref{eq:4})-(\ref{eq:6}); $\mathbf{\zeta(\mathbf{\textit{Q}_\textit{j}})} \gets \mathbf{\textit{C}_\textit{j}}$, \textit{Objective function} via Eq. (\ref{eq:35a})

\STATE {\textsc{Direction finding}}: $\mathbf{\textit{C}_\textit{j}^*} \gets \mathbf{\textit{C}_\textit{j}} $, via Eq. (\ref{eq:32});
\STATE \hspace{0.5cm}{$\mathbf{\textit{V}_\textit{j}} \gets \mathbf{\textit{C}_\textit{j}^*}$}, shifting all users from the bad to the best option choices \cite{ref33}

\STATE {\textsc{Step-size determination}}: $\mathbf{\theta_\textit{j}^*} $, via Golden section search method; $\mathbf{\sigma_\textit{j}} \gets \mathbf{\theta_\textit{j}^*}$ via Eqs. (\ref{eq:36})-(\ref{eq:38})
\STATE {\textsc{Update}}: $\mathbf{\textit{Q}_\textit{j+1}} = \mathbf{\textit{Q}_\textit{j} + \mathbf{\sigma_\textit{j}}(\mathbf{\textit{V}_\textit{j}} - \mathbf{\textit{Q}_\textit{j}})}$; $\mathbf{\zeta(\mathbf{\textit{Q}_\textit{j+1}})} \gets \mathbf{\textit{Q}_\textit{j+1}}$, repeat the \textsc{Evaluation} step (Line 3)
\STATE {\textsc{Convergence test}}: 
    \IF{$\mathbf{\zeta(\mathbf{\textit{Q}_\textit{j+1}})} \leq \mathbf{\zeta(\mathbf{\textit{Q}_\textit{j}})}$} 
    \STATE $\mathbf{\mathbf{\textit{Q}_\textit{1}^\textit{*}}}=\mathbf{\textit{Q}_\textit{j+1}}$, 
        $j = j+1$, 
        $\mathbf{\textit{Q}_\textit{j}} = \mathbf{\textit{Q}_\textit{1}^*}$
    \ELSE 
        \STATE{break}
    \ENDIF

\ENDFOR
\RETURN $\mathbf{\mathbf{\textit{Q}_\textit{1}}^\textit{*}}$
\STATE \textbf{End Procedure} (1)

\STATE \textbf{Procedure} {(2) OD-based flow shifting loop: }({\textsc{Input}}: $\mathbf{\textit{Q}_\textit{1}^\textit{*}}$ from \textbf{Procedure} (1); \textsc{Output}: $ \mathbf{\textit{Q}^\textit{*}}, \mathbf{\zeta(\mathbf{\textit{Q}^\textit{*}})}$ )\\
\FOR {$i \gets 0$ to $N$}
\STATE {\textsc{Evaluation}}: $\mathbf{\textit{M}_\textit{i}} \gets \mathbf{\textit{Q}_\textit{i}} $; $\mathbf{\textit{A}_\textit{i}} \gets \mathbf{\textit{M}_\textit{i}}$; $\mathbf{\textit{C}_\textit{i}} \gets \mathbf{\textit{A}_\textit{i}}$; $\mathbf{\zeta(\mathbf{\textit{Q}_\textit{i}})} \gets \mathbf{\textit{C}_\textit{i}}$
\STATE {\textsc{select randomly:}} $k \in K$
\STATE {\textsc{Direction finding}}(\textit{only within OD pair $k$}): $\mathbf{\textit{C}_\textit{i}^*(\textit{k})} \gets \mathbf{\textit{C}_\textit{i}(\textit{k})} $, {$\mathbf{\textit{V}_\textit{i}(\textit{k})} \gets \mathbf{\textit{C}_\textit{i}^*(\textit{k})}$}
\STATE {\textsc{Step-size determination}}: (\textit{only within OD pair $k$}): $\mathbf{\theta_\textit{i}^*(\textit{k})}$; $\mathbf{\sigma_\textit{i}(\textit{k})} \gets \mathbf{\theta_\textit{i}^*(\textit{k})}$ via Eqs. (\ref{eq:36}), (\ref{eq:38}), (\ref{eq:42})
\STATE {\textsc{Update}}: $\mathbf{\textit{Q}_\textit{i+1}} = \mathbf{\textit{Q}_\textit{i} + \mathbf{\sigma_\textit{i}(\textit{k})}(\mathbf{\textit{V}_\textit{i}} - \mathbf{\textit{Q}_\textit{i}})}$; $\mathbf{\zeta(\mathbf{\textit{Q}_\textit{i+1}})} \gets \mathbf{\textit{Q}_\textit{i+1}}$, repeat the \textsc{Evaluation} step (Line 20)

\STATE {\textsc{Convergence test}}: 

    \IF{$\mathbf{\zeta(\mathbf{\textit{Q}_\textit{i+1}})} \leq \mathbf{\zeta(\mathbf{\textit{Q}_\textit{i}})}$ } 
    \STATE $\mathbf{\mathbf{\textit{Q}_\textit{2}^\textit{*}}}=\mathbf{\textit{Q}_\textit{i+1}}$, 
        $i = i+1$, 
        $\mathbf{\textit{Q}_\textit{i}} = \mathbf{\textit{Q}_\textit{2}^*}$
    \ELSIF{$i \leq N_2$}
    \STATE {$i = i+1$, $\mathbf{\textit{Q}_\textit{i}} = \mathbf{\textit{Q}_\textit{2}^*}$}
    \ELSE 
        \STATE{break}
    \ENDIF
\ENDFOR
\RETURN $ \mathbf{\textit{Q}_\textit{2}^\textit{*}}$, $\mathbf{\zeta(\mathbf{\textit{Q}_\textit{2}^\textit{*}})}$, $\mathbf{\textit{Q}^\textit{*}} = \mathbf{\textit{Q}_\textit{2}^*}$, $\mathbf{\zeta(\mathbf{\textit{Q}^\textit{*}})} = \mathbf{\zeta(\mathbf{\textit{Q}_\textit{2}^\textit{*}})}$
\STATE \textbf{End Procedure} (2)
\end{algorithmic}
\end{algorithm}

The system-based flow shifting strategy uses an `all-at-once' approach, simultaneously updating flows for all OD pairs with the system step size and OD relative gaps $\mathbf{\partial} = \{\partial_k, \forall k \in K\}$, and time option ratios $\varphi = (\varphi_{k}(t),  \forall k \in K, t \in \mathbf{\tau}_k)$, where $\mathbf{\tau}_k$ is the set of non-optimal time options. The OD-based flow shifting strategy employs an `one-at-a-time' approach, updating flows for a single OD pair per iteration. The AdaGDD algorithm is designed to find a local optimal solution for the objective function in Eq. (\ref{eq:35a}), and this local optimal solution is then used as the DTUE-PT solution.

From the perspective of system-based flow shifting loop, it simultaneously updates solutions for all OD pairs using the system step size $\theta$, $\mathbf{\partial}$, and $\boldsymbol\varphi$ to update the actual step size $\mathbf{\sigma} = \{\sigma_{k}(t), \forall k \in K, t \in T_{o(k)}\}$. These $\sigma_{k}(t)$, $\partial_k$ and $\varphi_{k}(t)$ are updated by Eqs. (\ref{eq:36})-(\ref{eq:38}).
\begin{equation}
\label{eq:36}
\sigma_{k}(t) = \theta * \partial_k * \varphi_{k}(t), \quad \forall k \in K, t \in T_{o(k)}
\end{equation}

Eq. (\ref{eq:36}) is used to calculate the actual step size $\sigma_{k}(t)$ (which determines the number of users to be shifted for option $(k, t)$). This calculation involves the step size $\theta$ (from the system perspective to determine the optimal system step size), OD relative gap $\partial_k$ (the gap relative to the `ideal' UE for option $(k, t)$, where all used time options have the same and lowest option cost), and the time option ratio $\varphi_{k}(t)$ (to determine how to target users in non-optimal options for shifting).

\begin{equation}
\label{eq:37}
\partial_k = \frac{\displaystyle C_k(\tau_k)- C_k^*}{C_k(\tau_k)}, \quad \forall k \in K
\end{equation}

In Eq. (\ref{eq:37}), $C_k^*$ is the minimum option cost within OD pair $k$, which is calculated by Eq. (\ref{eq:32}). This is used to measure the gap between the average cost and the optimal cost within OD pair $k$. For any given value of $C_k^*$, $\partial_k$ is minimized by minimizing $C_k(\tau_k)$. A larger gap will stimulate a greater proportion of users to be shifted within this OD pair $k$. 

\begin{equation}
\label{eq:39}
C_k(\tau_k) = \frac{ \displaystyle \sum_{t \in T_{o(k)}}C_k(t)}{n}, \quad \forall k \in K
\end{equation}
where $C_k(\tau_k)$ is the average time option cost of all non-optimal time options. $n$ represents the total number of available time options within OD pair $k$.

\begin{equation}
\label{eq:38}
\varphi_{k}(t) = \frac{C_k(t)}{{\displaystyle \sum_{t_1 \in T_{o(k)}}C_k(t_1)}}, \quad \forall k \in K, t \in \mathbf{\tau}_k
\end{equation}

Eq. (\ref{eq:38}) indicates that the number of users from non-optimal options waiting to be shifted depends on the weighted ratio of the time options within OD pair $k$. A non-optimal time option with a higher option cost will result in a greater number of users being shifted to the optimal option. For example, there is public transport network with 2 OD pairs:  In OD pair OD1, there are 3 time options (option 1-1 (option cost 10), option 1-2 (option cost 20), option 1-3 (option cost 15)); In OD pair OD2, there are another 3 time options (option 2-1 (option cost 50), option 2-2 (option cost 60), option 2-3 (option cost 40)). Then OD relative gap and time option ratios are $\mathbf{\partial} = \{ 0.33, 0.20\}$, which is calculated by (OD1: $(15-10)/15 = 0.33$, OD2: $(50-40)/50 = 0.20$), and $\mathbf{\varphi} =\{(0.57, 0.43), (0.45, 0.55)\}$, which is calculated by (Option 1-2: $20/(20+15) = 0.57$, Option 1-3: $15/(20+15) = 0.43$), and (Option 2-1: $50/(50+60) = 0.45$, Option 2-2: $60/(50+60) = 0.55$). 

From the perspective of the OD-based flow shifting loop, its solution is updated for a single OD pair per iteration. Most steps of this loop are the same as those in the system-based flow shifting loop, except for the determination of the step size. Eq. (\ref{eq:37}) is replaced with:\\
\begin{equation}
\partial_k =  \begin{cases} 
    1, \quad \text{if } k \text{ is targeted}\\
    \\
    0, \quad  \text{if } k \text{ is not targeted}
\end{cases}
\label{eq:42}
\end{equation}\\
The OD-based flow shifting loop ends when transferring just one user from a non-optimal option to any optimal option, on any OD, fails to improve the system gap. 


\section{Case study}\label{sec:case-study}
This study proposes the AdaGDD solution algorithm to address the NMP problem and validates its performance by comparing it with the MSA method \cite{ref34} and the day-to-day learning method \cite{ref27}. 

In this section, we first conduct experiments on a small network to examine the properties of the AdaGDD solution algorithm, as detailed in Section \ref{sec:smallnetwork}. The experiments in Section \ref{sec:result1} is designed to illustrate the convergence properties of these three methods, while Section \ref{sec:result2} compares their sensitivity by testing different initial model solution settings. Finally, we apply the proposed method to the extracted Hong Kong Mass Transit Railway (MTR) network to explore its application in network design in Section \ref{sec:result3} and discuss its extension to incorporate route choice in Section \ref{sec:result4}.

\subsection{Synthetic network}\label{sec:smallnetwork}
\begin{figure}[ht]
    \centering
    \includegraphics[width=0.496\textwidth]{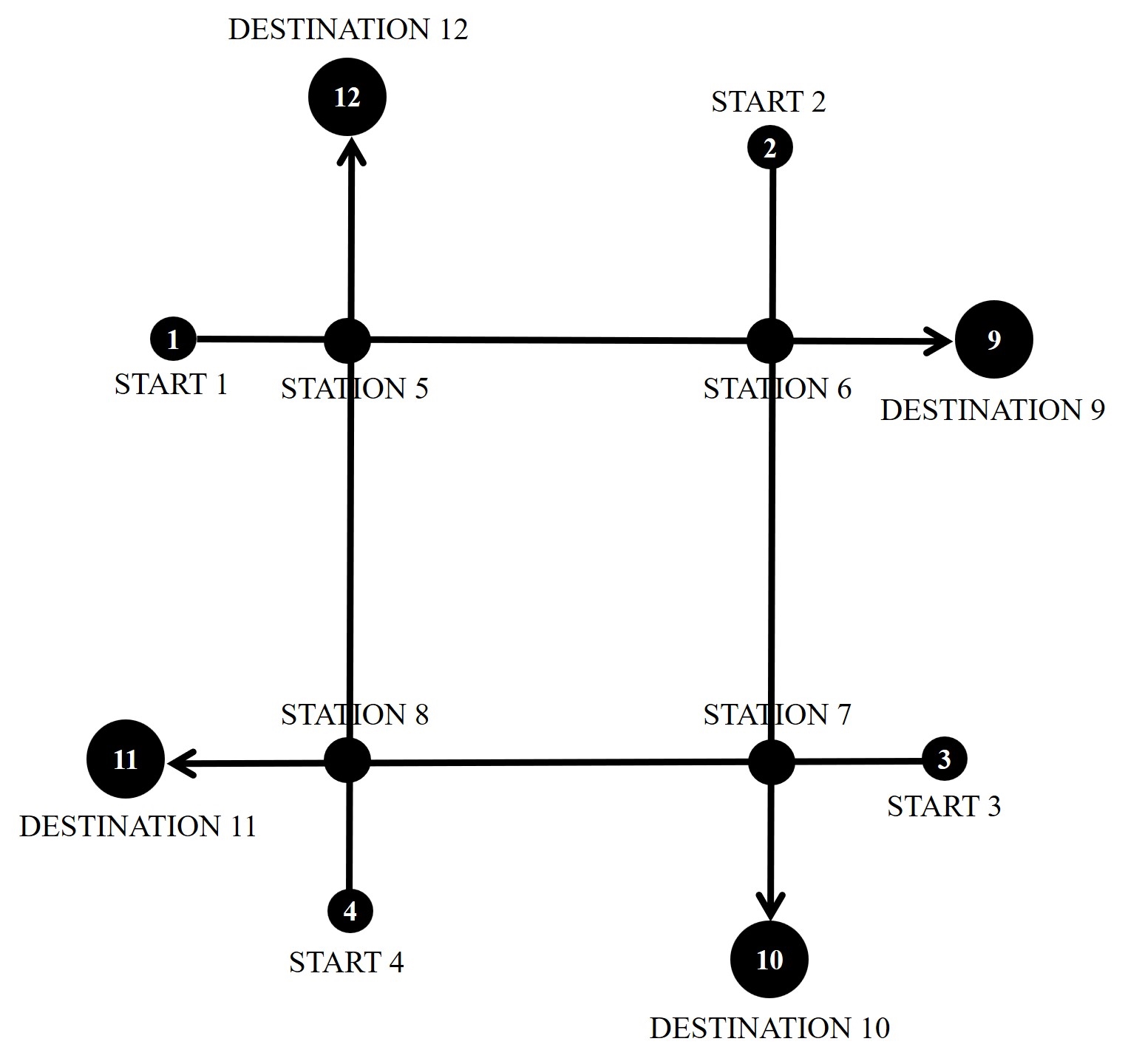}
    \caption{Test Public Transport Network}
    \label{fig:fig5}
\end{figure}

We validate the method performance using a multi-line network with transfers, and crowding at key stations, and comparing with commonly used day-to-day learning and MSA methods. Fig. \ref{fig:fig5} shows the public transport network. It comprises 4 lines traversing various stations: \textit{Line 1} (1-5-6-9), \textit{Line 2} (2-6-7-10), \textit{Line 3} (3-7-8-11), and \textit{Line 4} (4-8-5-12). Each line operates with a 5-minute headway and capacity of 230 users per train from 5:00 AM onward. 

\begin{figure*}[ht]
    \centering
    \includegraphics[width=0.92\textwidth]{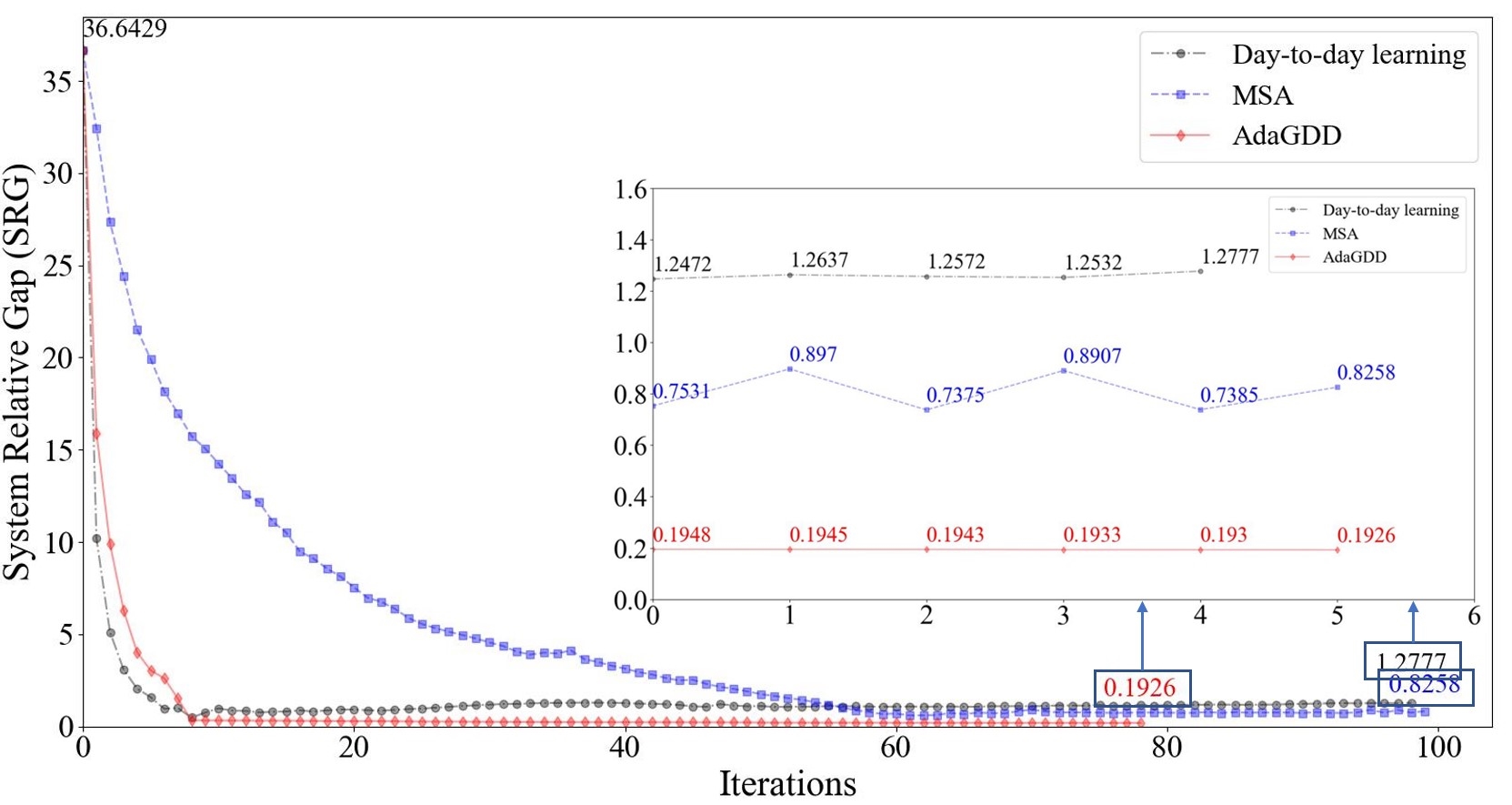}
    \caption{Comparison of model solution algorithms: Day-to-day learning, MSA, and AdaGDD}
    \label{fig:fig6}
\end{figure*}
We consider 16 OD travels for 4 origins (1, 2, 3, and 4) and 4 destinations (9, 10, 11, and 12). Take the OD pair `1-11' for example: the user's route is: \textbf{1} (line 1) - \textbf{6} (transfer from line 1 to line 2) - \textbf{7} (transfer from line 2 to line 3) -\textbf{11} (line 3). We set the OD demand ($K$) to be 2000 users per OD pair with a desired arrival time at the destination (e.g., commence work) at 9:00 AM.  Table \ref{tab:table-6} summarizes the preferred departure times for different OD pairs.
\begin{table}[htbp]
    \begin{minipage}[b]{0.498\textwidth}
      \centering
      \caption{Initial Preferred departure times (AM) for ODs}
      \label{tab:merged-cells-academic}
      \begin{tabular}{cccccc}
        \toprule
        \textbf{Origins} & \multicolumn{3}{c}{\textbf{Destinations}} \\
        \cmidrule(lr){2-5}
        &  \textbf{9} & \textbf{10} & \textbf{11} & \textbf{12} \\
        \midrule
        \textbf{1} & 8:50:00 & 8:35:00 & 8:30:00 & 8:45:00\\
        \textbf{2} & 8:45:00 & 8:35:00 & 8:30:00 & 8:20:00 \\
        \textbf{3} & 8:25:00 & 8:40:00 & 8:40:00 & 8:30:00\\
        \textbf{4} & 8:30:00 & 8:15:00 & 8:40:00 & 8:35:00 \\
       
        \bottomrule
      \end{tabular}
      \label{tab:table-6}
    \end{minipage}
\end{table}

This small but complex network is designed to serve the purpose of obtaining validation requirements while maintaining parsimony. It includes all key service network features in real-world: multiple OD demands, multiple service lines, multiple transfers across lines (ranging from 0 to 3 changes), and denied boarding at boarding stations. 
For the model setting, we set the waiting, early arrival, and late arrival penalties to be 10, 1, and 10 respectively. To evaluate the performance of the solution algorithms, the System Relative Gap (SRG) is employed. A smaller SRG indicates better performance, as outlined in \cite{ref24}:\\
\[ SRG = \frac{\displaystyle \sum_{k \in K} \sum_{t \in T_{o(k)}} (C_k(t) - C_k^*) \times q_k(t)}{ \displaystyle \sum_{k \in K} \sum_{t \in T_{o(k)}} q_k(t) \times C_k^*}\]\\
where the denominator calculates the `ideal' system travel cost, assuming all users experience the least travel cost, $C_k^* = \min_{t \in T_{o(k)}}C_k(t)$. The numerator calculates the absolute gap between the model solution cost to the `ideal' system travel cost. In essence, the SRG measures the level of the deviation of the model solution from the `ideal' solution. The SRG facilitates the analysis and comparison of performance among different solution models. A smaller SRG indicates a better solution (closer to the `ideal' UE solution), as average travel costs of different departure time options are closer to each other.

\subsubsection{Model validation}\label{sec:result1}\

Fig. \ref{fig:fig6} shows the convergence performance of compared algorithms, day-to-day learning, MSA and AdaGDD methods, using the default initial model solution setting (see Table \ref{tab:initial-settings}). The AdaGDD and day-to-day learning algorithms converge quickly (around 10 iterations), while the MSA needs more iterations to converge. Furthermore, the AdaGDD algorithm converges to a better solution with the smallest relative gap (SRG) of 0.1926, which improves the benchmark solution by 85\% for day-to-day learning (1.2777) and 76\% for the MSA (0.8258). 

\subsubsection{Sensitivity analysis}\label{sec:result2}\
\begin{figure*}[ht]
    \centering
    \includegraphics[width=0.81\textwidth]{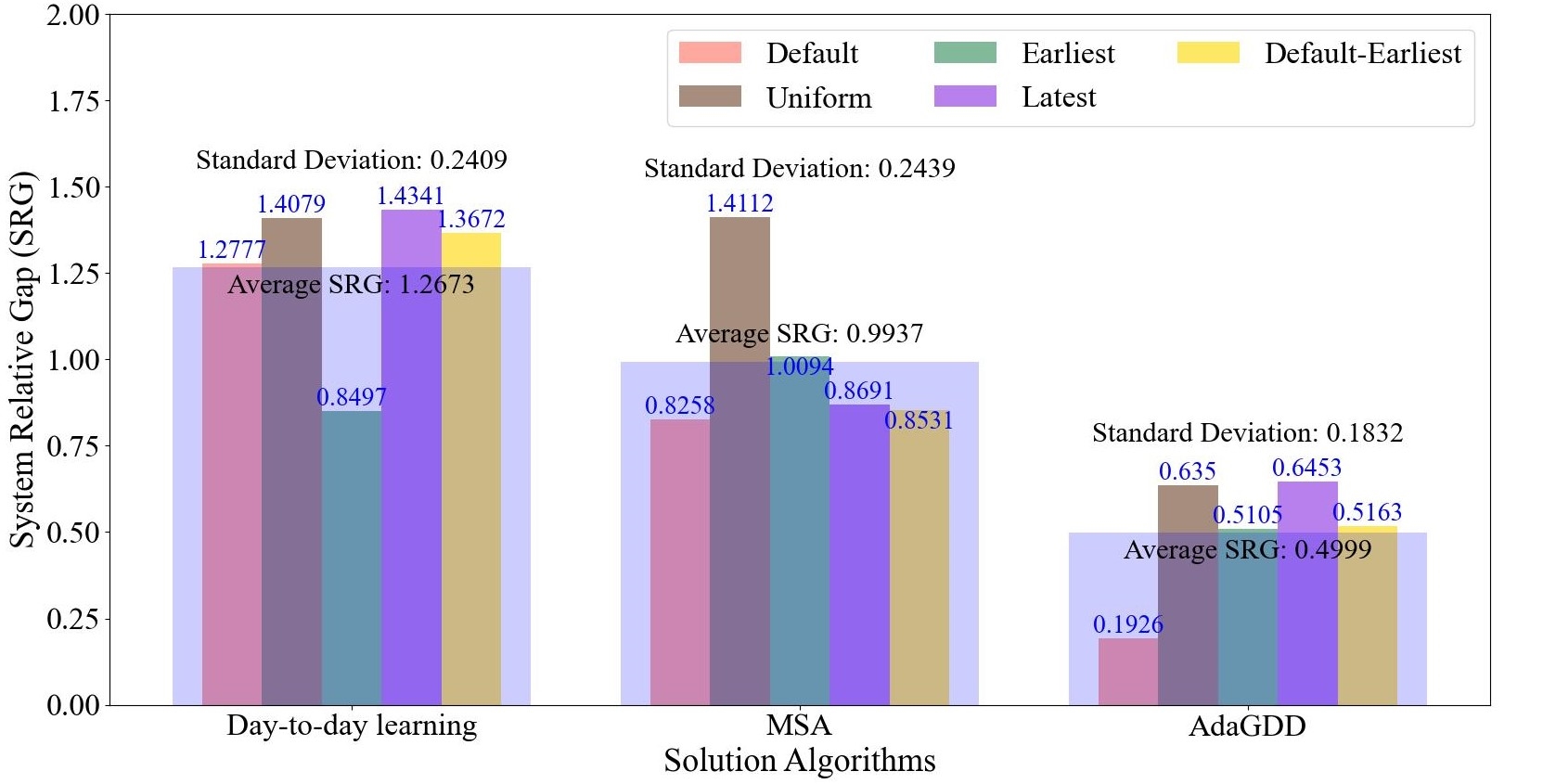}
    \caption{System relative gaps (SRG) under 5 initial model solution settings: Day-to-day learning , MSA, and AdaGDD}
    \label{fig:fig14}
\end{figure*}

We tested five different initial model solution settings using day-to-day learning, MSA, and AdaGDD methods respectively. It covers a broad range of initial model solution settings including some extreme condition (users choose the latest train i.e., 12:25 PM which is impossible in reality owing to unbearable lateness to work starting time 9:00 AM) to further test the algorithm robustness performance. 
\begin{itemize}
    \item \textit{Default:} All users depart at the preferred times (see Table \ref{tab:table-6}) assuming infinite service capacity;
    \item \textit{Uniform:} All users are distributed evenly at all time options considered in the study; 
    \item \textit{Earliest:} All users depart at the earliest time option (i.e., 5:00 AM); 
    \item \textit{Latest:} All users depart at the latest time option (i.e., 12:25 PM); 
    \item \textit{Default-Earliest:} Half users depart at the default times and half at the earliest time option;
\end{itemize}

Detailed information of these 5 initial model solution settings can be found in Table \ref{tab:initial-settings} (Appendix-D).

Fig. \ref{fig:fig14} shows the system relative gaps (SRG) 
under five initial solution settings using these three solution algorithms. We use the average SRG to represent the SRG of that algorithm, which is filled with light color over five results in Fig. \ref{fig:fig14}. It shows that the AdaGDD can achieve an average SRG of 0.4999, improving the benchmark performance by 60.6\% for day-to-day learning (1.2673) and 49.7\% for MSA (0.9937). This further supports the findings in previous section that the AdaGDD can achieve the best DTUE-PT solution, which improves the benchmark performances (day-to-day learning and MSA). The result also shows that these 3 algorithms perform well/poorly in some specific initial model solution setting(s) (which is/are line(s) far away from the average SRG). The day-to-day learning method can find its best solution with SRG 0.8497 for initial model solution setting `Earliest'. The MSA algorithm performs poorly with a solution of SRG 1.4112 for initial model solution setting `Uniform'. The AdaGDD method performs well with a solution of SRG 0.1926 for initial model solution setting `Default'. From the perspective of best solutions, the AdaGDD has great potential to find the best DTUE-PT solution with SRG 0.1926, which is the lowest, compared with both best solutions from day-to-day learning (SRG 0.8497) and MSA (SRG 0.8258). From the perspective of the worst solutions, the AdaGDD only shows the worst DTUE-PT solution with SRG 0.6453, which is still the lowest, compared with both worst solutions from day-to-day learning (SRG 1.4341) and MSA (SRG 1.4112). From the perspective of the standard deviation, the AdaGDD shows the lowest value 0.1832, which is the least compared with 0.2409 (day-to-day) and 0.2439 (MSA). Therefore, we can conclude that the proposed AdaGDD can achieve the best DTUE-PT solution with the least SRG and show the best robustness performance with the least standard deviation 0.1832.

\subsection{Hong Kong MTR network} \label{sec:result3}

\begin{figure}[ht]
    \centering
    \includegraphics[width=0.49\textwidth]{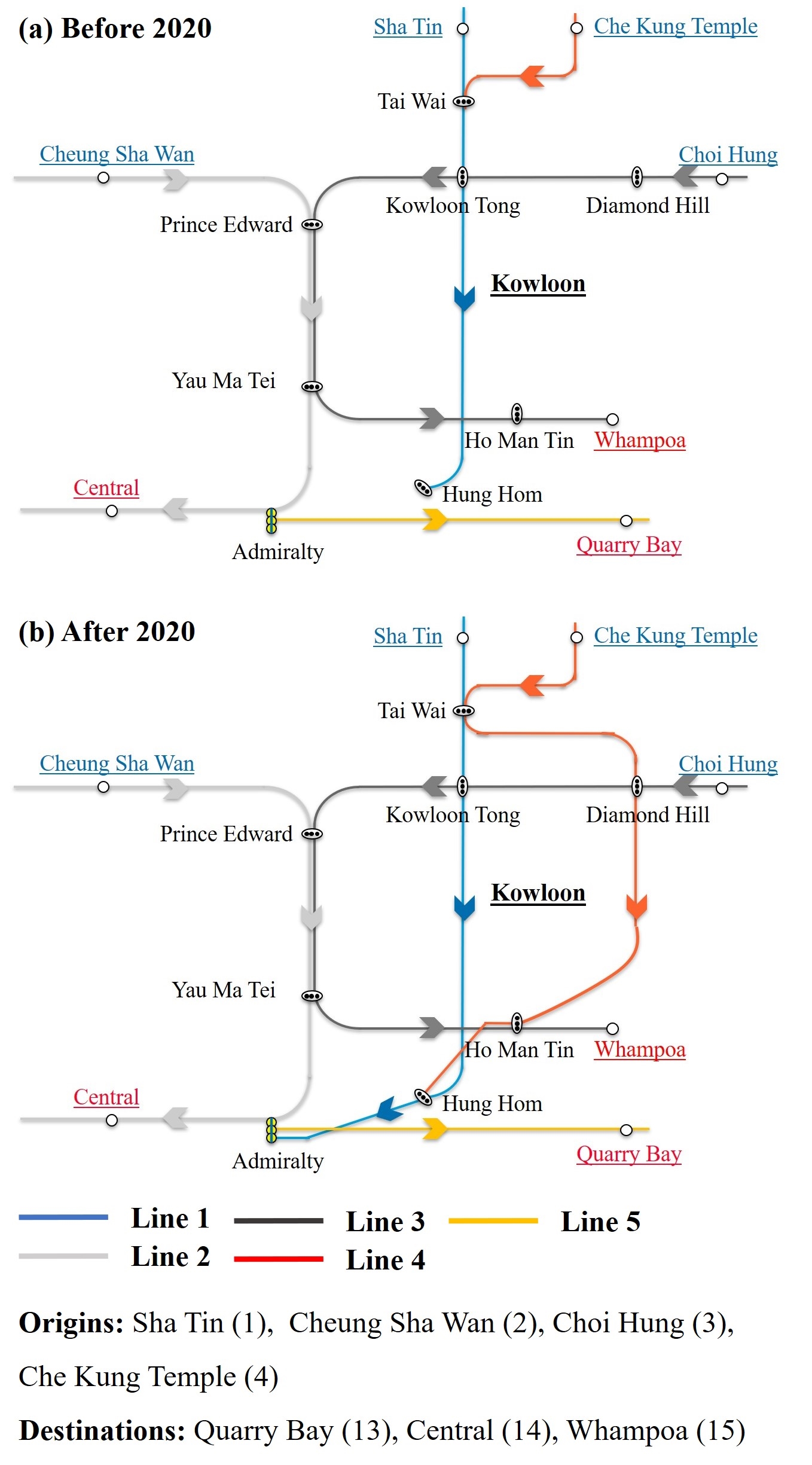}
    \caption{Extracted Hong Kong transit networks: (a) previous network (before 2020), (b) current network. 
    Reproduced from web images (https://www.checkerboardhill.com/2020/05/mtr-rail-operating-divisions/) and (https://www.mtr.com.hk/en/customer/services/system\_map.html)}
    \label{fig:24}
\end{figure}

We further implement the AdaGDD method in the extracted Hong Kong MTR network to examine the effects of network design on system performance, as shown in Fig. \ref{fig:24}. In this paper, we model the central part of the Hong Kong MTR network. We consider specific OD demand flows and lines running from north to south, towards the CBD. Other demand flows within the CBD are not included in this paper due to the unavailability of data. If the data becomes available, these demand flows can be incorporated as temporal background flow to update the corresponding train capacities in our model.

Fig. \ref{fig:24}(a) shows the network before 2020. The new stretch of the Tuen Ma Line (red line) from Tai Wai to Admiralty was completed on June 27, 2021, and the East Rail Line (blue line) was extended from Hung Hom in Kowloon to Admiralty on Hong Kong Island on May 15, 2022, as shown in Fig. \ref{fig:24}(b). 
The modeled section includes 5 lines indicating 26 trains (line 1), 26 trains (line 2), 16 trains (line 3), 35 trains (line 4), and 17 trains (line 5) respectively. Each train has a capacity $S_t =2600$ users. The OD matrix (see Table \ref{tab:OD}) is derived from card data, comprising a total of 52,717 trips. Real train schedules from Google Maps are used in this network. The weighting parameters used in Eqs. (\ref{eq:4})-(\ref{eq:6}) are $\alpha = 18, \beta = 5$, and $\gamma = 12$ \cite{ref35}.

\begin{table}[ht]
    \begin{minipage}[b]{0.49\textwidth}
        \centering
        \caption{Origin-Destination Demand from smart card (5.50-10.00 AM)}
        \label{tab:OD}
        \begin{tabularx}{\columnwidth}{cX X X}
            \toprule
            \textbf{O/D} & \textbf{Quarry 13} & \textbf{Central 14} & \textbf{Whampoa15} \\
            \midrule
            \textbf{Sha Tin-1} & 5356 & 5663  & 1892  \\
            \textbf{Cheung Sha Wan-2} & 6116  & 4073 & 2049  \\
            \textbf{Choi Hung-3} & 14967 & 5852  & 2525  \\
            \textbf{Che Kung Temple-4} & 1727 & 1848  & 649 \\
            \bottomrule
        \end{tabularx}
    \end{minipage}
\end{table}
We plan to measure the system performance of these two network designs at different demand levels to identify the potential outcomes of the new network design (current network). Demand level refers to the relative demand compared to the real demand (demand level 100\%) in Table \ref{tab:OD}. For example, a 165\% demand level means multiplying the real demand by 1.65, and a 60\% demand level means multiplying the real demand by 0.6. The 165\% demand level corresponds to a scenario where future city development may stimulate increased travel demand. The 60\% demand level corresponds to a scenario where, in some cities, newly built districts attract demand from old districts, resulting in a lower demand level in these older areas.
\begin{figure}[ht]
    \centering
    \includegraphics[width=0.5\textwidth]{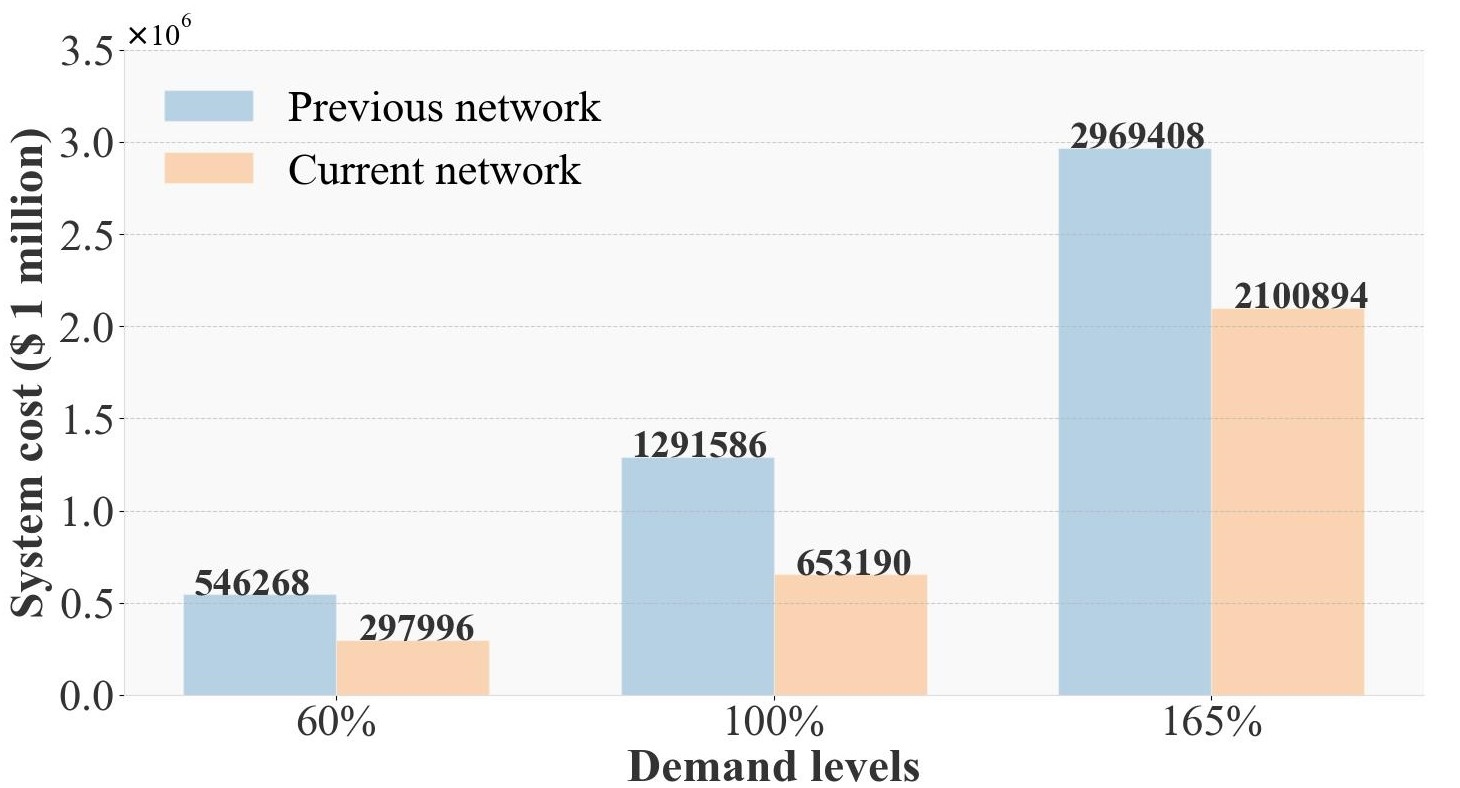}
    \caption{System costs of these two network designs: considering demand levels of 60\%, 100\%, 165\% using the AdaGDD method}
    \label{fig:25}
\end{figure}

Fig. \ref{fig:25} shows that the current network achieves better system performance with lower system costs compared to the previous network. The system cost improvement ratios (1 - Current/Previous) are 45.45\%, 49.43\%, and 29.25\% at demand levels of 60\%, 100\%, and 165\%, respectively. This indicates that the current network design performs well under low demand conditions (60\% and 100\% demand levels), but may be less efficient under high demand conditions (165\% demand level). This may be due to higher demand occupying most of the train capacity, which limits available choices for users and reduces the potential for system improvement. 

\begin{figure}[ht]
    \centering
    \includegraphics[width=0.5\textwidth]{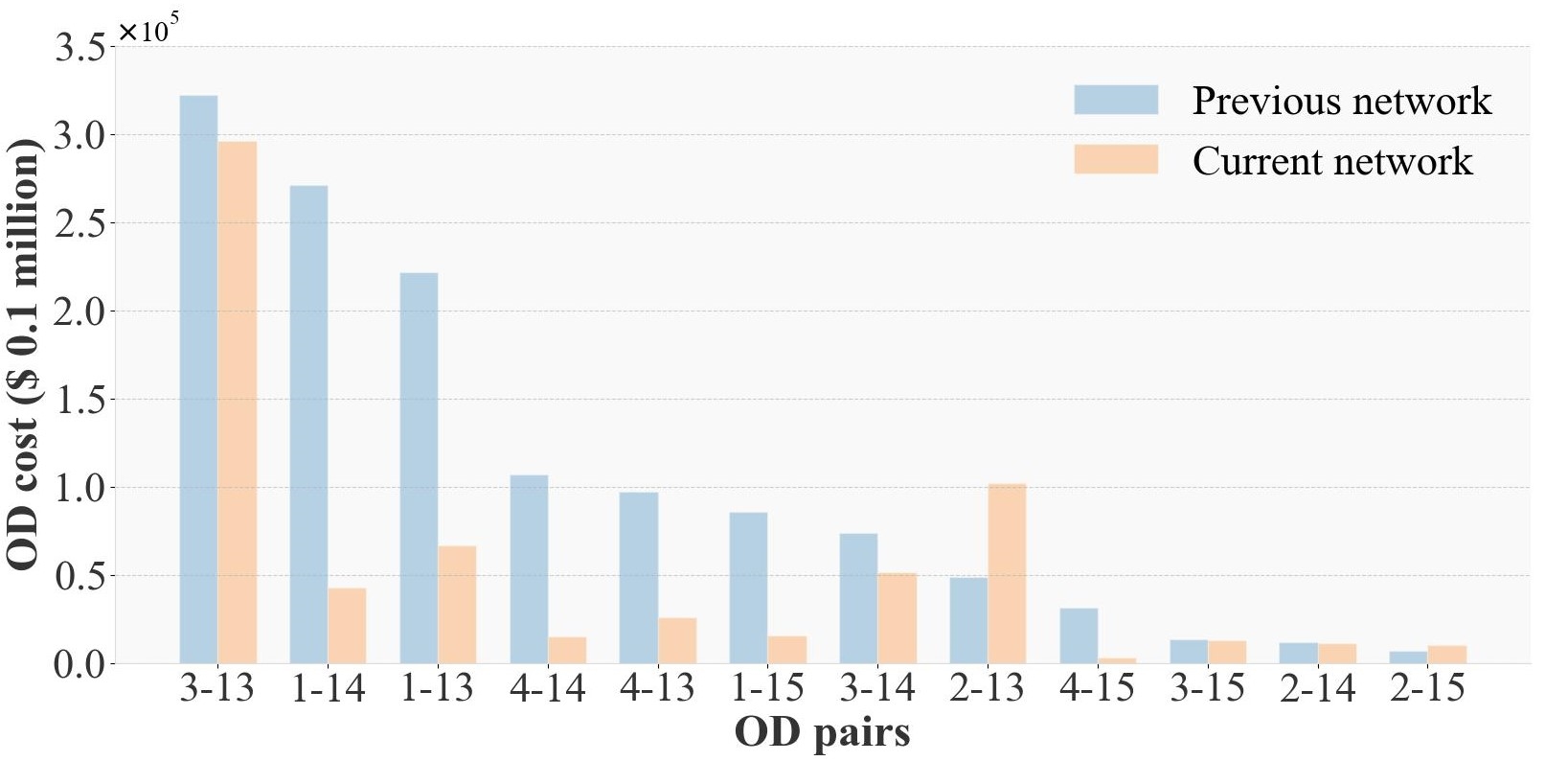}
    \caption{OD pair costs of these two network designs at the demand level 100\% using the AdaGDD method}
    \label{fig:26}
\end{figure}

To understand what happens among OD pairs in these two network designs, we take the 100\% demand level as an example. Fig. \ref{fig:26} shows the OD costs for the two network designs. The OD pair numbers correspond to the station numbers shown in Table \ref{tab:OD}
(e.g., OD `1-13' represents the OD pair `Sha Tin - Quarry Bay').
\begin{figure*}[ht]
    \centering
    \includegraphics[width=0.9\textwidth]{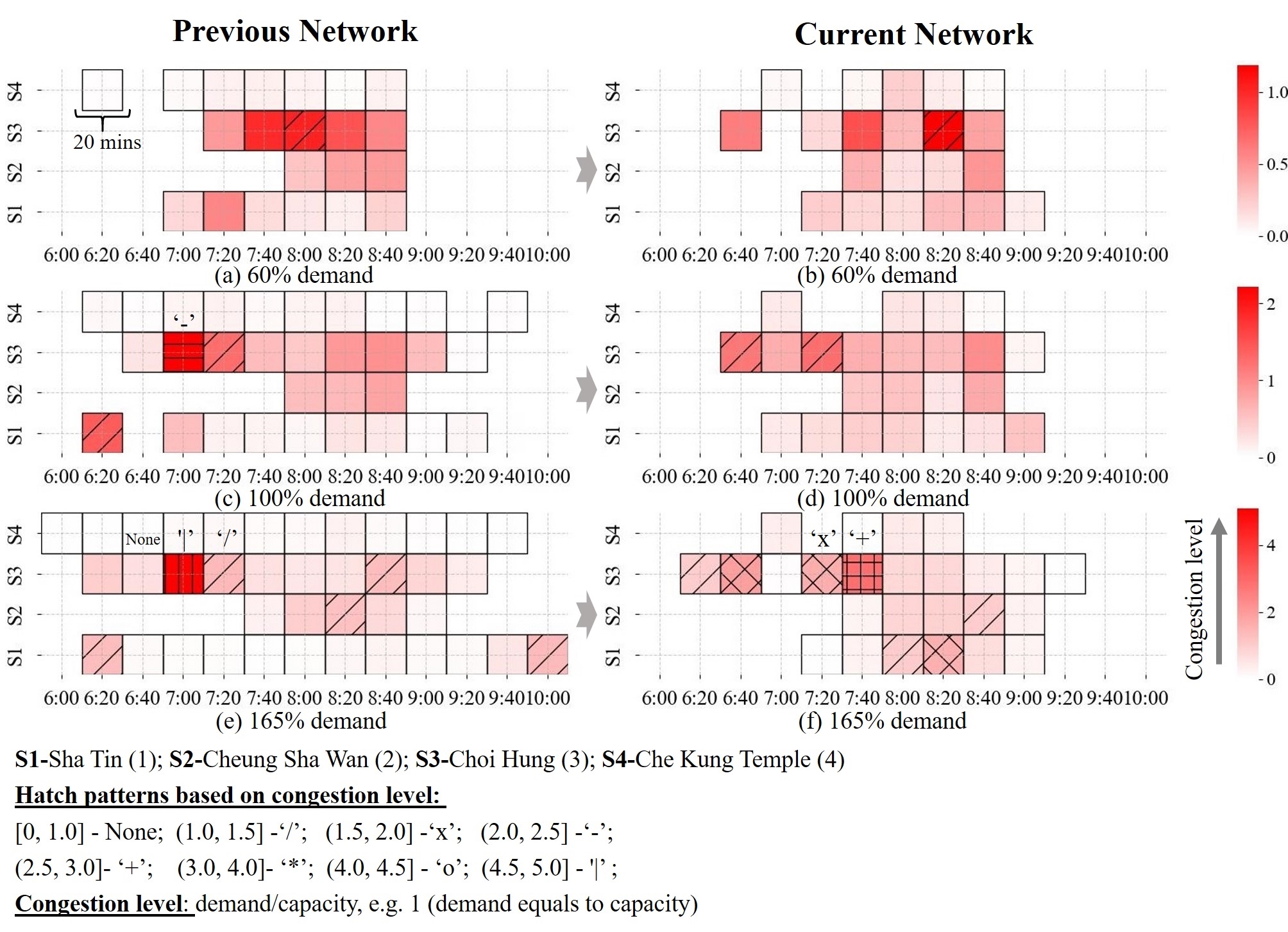}
    \caption{Congestion heat map of these two network designs: considering demand levels of 60\%, 100\%, 165\% using the AdaGDD method}
    \label{fig:27}
\end{figure*}
The result shows that most OD costs decrease when updating the network from the previous version to the current one, except for the OD pairs `2-13' and `2-15'. The decreased OD costs are as expected because some loads on the link `Prince Edward - Admiralty' are shifted to the new links `Diamond Hill - Hung Hom' (via Line 4) and `Kowloon Tong - Admiralty' (via Line 1). This will alleviate congestion on the critical link `Prince Edward-Admiralty'. The increased OD costs for the pairs `2-13' and `2-15' may be due to transfer requirements at stations `Admiralty' or `Prince Edward'. Users from `2-13' and `2-15' may experience higher costs because they are less competitive compared to upstream users from the origin stations `Che Kung Temple' and `Choi Hung' in the current network design. Users from OD pairs `4-15', `4-14', `4-13', `1-14', and `1-15' benefit the most from the current network design.

To explore the effects of network design on peak-hour congestion, we plot a heat map of the congestion state at the origin stations (Sha Tin, Cheung Sha Wan, Choi Hung, and Che Kung Temple) for each time period, as shown in Fig. \ref{fig:27}. The time interval is set to 20 minutes, starting from 5:50 AM to 10:10 AM. The congestion level is calculated by dividing the total demand by the total available train capacity for that time slot (e.g., a congestion level of 2 indicates that demand is twice the train capacity).

The results show that obvious congestion occurs at Choi Hung at a demand level of 60\%, at both Choi Hung and Sha Tin at a demand level of 100\% (in the previous network), and simultaneously at Choi Hung, Cheung Sha Wan, and Sha Tin at a demand level of 165\% respectively. The worst congestion occurs at Choi Hung station, where demand is significantly higher compared to other stations.  From the perspective of network design, in the current network, peak hours (the period with the highest congestion level) predominantly occur later compared to the previous network. For example, at a demand level of 165\%, the peak hour at Choi Hung station shifts from 7:00 AM to 7:40 AM. Overall congestion is also found to be alleviated (decreased congestion levels) in the current network, compared to the previous network. For example, at a demand level of 165\%, the highest congestion level decreases from 5 to 3. Furthermore, some peak-hour demands shift to neighboring time periods in the current network compared to the previous network. For example, at a demand level of 100\% at Choi Hung station, the peak-hour demand shifts from 7:00 AM in the previous network to 6:40 AM and 7:20 AM in the current network. This further supports the observed phenomenon of relieved congestion in the current network. In conclusion, the current network can alleviate congestion compared to the previous network, and the proposed method can help analyze the potential benefits of future network designs.

\subsection{Model Extension (incorporating route choices)} \label{sec:result4}
\begin{figure*}[ht]
    \centering
    \includegraphics[width=0.99\textwidth]{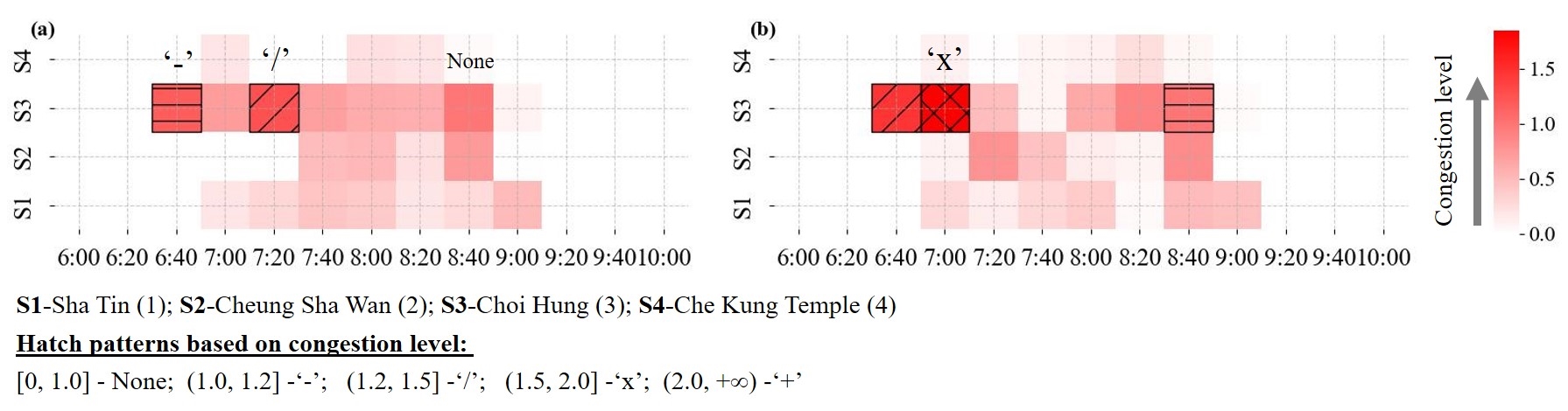}
    \caption{Congestion heat map of these two scenarios:  (a) Departure time only (left), and (b) Simultaneous Departure time and Route Choice (right)}
    \label{fig:30}
\end{figure*}
In this section, we extend our model further to incorporate both simultaneous departure time and route choices.  This can be achieved by simply extending the variables to include route choice. For example, $q_k(t)$ is updated to $q_k(t, r)$, where $r$ represents the route, $r \in R_k$ ($R_k$ is the set of routes for OD pari $k$) .  The framework for achieving the UE solution remains exactly the same as described above. The detailed updated model is provided in Appendix-B.
\begin{table}[ht]
    \centering
    \begin{minipage}[t]{0.46\textwidth}
        \centering
        \caption{Origin-Destination Matrix(5:50-10:00 AM)}
        \label{tab:top_table}
        \begin{tabularx}{\columnwidth}{cX X X}
            \toprule
            \textbf{O/D} & \textbf{Quarry Bay} & \textbf{Central} & \textbf{Whampoa} \\
            \midrule
            Sha Tin & 5356 (S) & 5663 (D) & 1892 (D) \\
            Cheung Sha Wan & 6116 (S) & 4073 (S) & 2049 (S) \\
            Choi Hung & 14967 (D) & 5852 (S) & 2525 (S) \\
            Che Kung Temple & 1727 (D) & 1848 (D) & 649 (D) \\
            \bottomrule
        \end{tabularx}
    \end{minipage}
\end{table}

We test the same transit system (Fig. \ref{fig:24}(b) - current network design) as described in Section \ref{sec:result3}. For certain OD pairs, a second route (`D' in Table \ref{tab:top_table}) is provided if two routes have similar transfer requirements for that OD pair. Otherwise, only a single route (`S' in Table \ref{tab:top_table}) is available.

As shown in Fig. \ref{fig:30}, compared to the scenario with only departure time, incorporating route choice leads to passenger accumulation during certain time periods (e.g., 6:30 AM–7:30 AM at Choi Hung, 7:10 AM–7:30 AM at Cheung Sha), rather than an even distribution across neighboring time periods. Therefore, providing more choices to passengers may lead to increased congestion if no additional interventions (e.g., incentives) are implemented.


\begin{figure}[ht]
    \centering
    \includegraphics[width=0.48\textwidth]{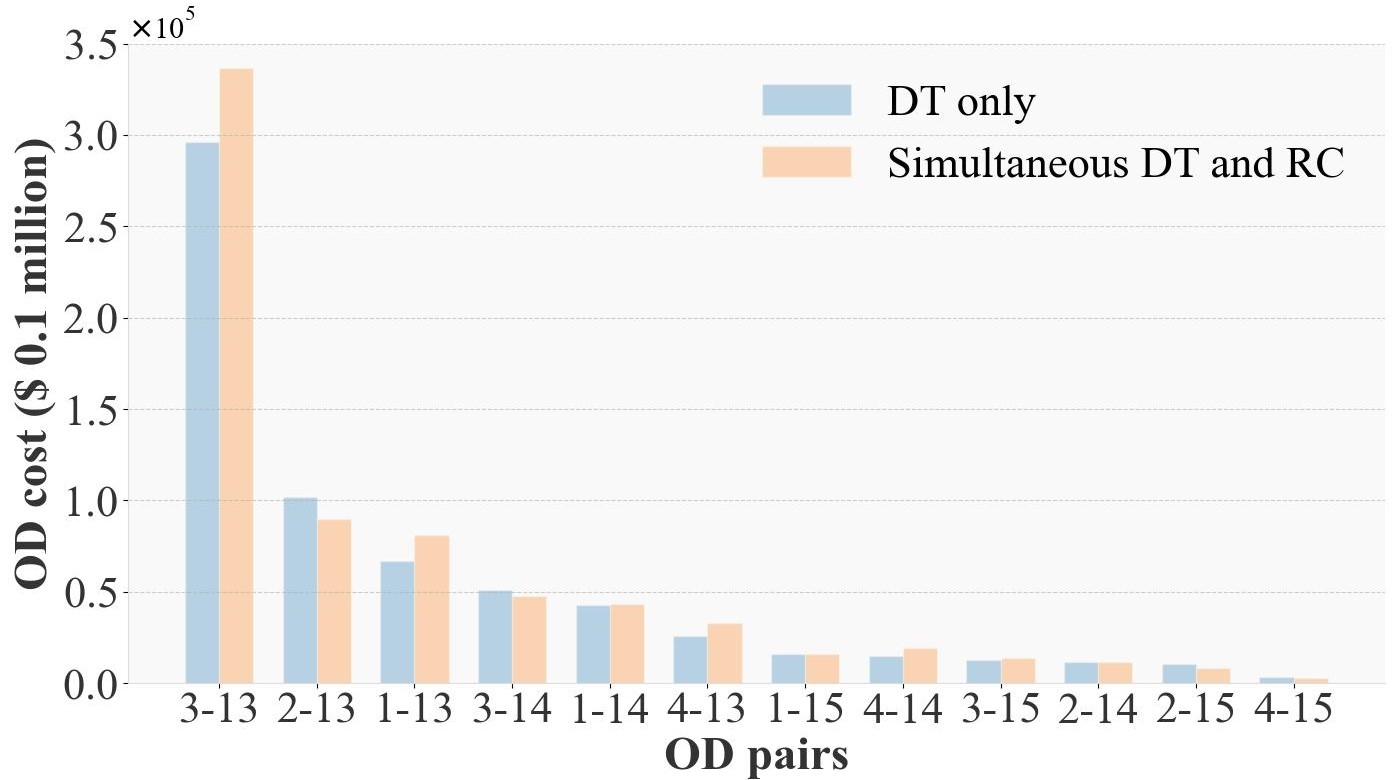}
    \caption{OD pair costs of these two scenarios: Departure time only and Simultaneous Departure time and Route Choice}
    \label{fig:29}
\end{figure}

Fig. \ref{fig:29} illustrates that providing additional route choices for users worsens the congestion state for most OD pairs (e.g., `3-13', `1-13', `1-14', '4-13', `4-14', `3-15'), resulting in increased OD costs. This outcome is expected, as passengers act selfishly and make decisions based on their individual benefits under the UE framework. When more flexibility is provided, passengers are likely to make fewer `good' decisions as anticipated by the system, leading to greater inefficiencies. This finding may suggest a situation analogous to the Braess paradox as previously discussed in \cite{ref69}.

\section{Conclusion} \label{sec:conclusion}
This paper introduces the departure time choice user equilibrium problem in public transport (DTUE-PT), assuming users can not reduce their perceived travel costs by changing departure times unilaterally. This paper considers the DTUE-PT problem for multi-line, schedule-based networks with hard train capacity constraints. The DTUE-PT problem is formulated as a Non-linear Mathematical Program (NMP) model to estimate the demand flow distribution over departure times by minimizing the system gap of travel cost between non-optimal options and the optimal option, with the system dynamics and passenger interactions modeled using a simulation. We develop an adaptive gap-based descent direction (AdaGDD) algorithm to solve the NMP model. We validate the proposed methodology performance by comparing with traditional Method of Successive Average (MSA) and day-to-day learning methods in public transport assignment problems and test its sensitivity and robustness performance by testing different initial model solution settings. We also apply the proposed method to evaluate its potential application in network design using the Hong Kong MTR network. The DTUE-PT model can be easily extended to include route choice and solved using the same framework. 

The result shows that our proposed methodology has great potential to find the accurate and robust DTUE-PT solution. The AdaGDD achieves a better DTUE-PT solution with the least system relative gap 0.1926, demonstrating a significant improvement compared to day-to-day learning and MSA, yielding 85\% and 76\% better results, respectively. The sensitivity analysis highlights that the AdaGDD method outperforms the day-to-day learning and MSA by providing a minimal standard deviation of 0.1832, compared to 0.2409 (day-to-day learning) and 0.2439 (MSA), respectively. The proposed method demonstrates great potential for evaluating system performance in new network designs and reveals that the current Hong Kong MTR network does alleviate congestion compared to the previous network. Providing greater flexibility to passengers has been found to worsen system congestion in the UE state.

This DTUE-PT study presents a valuable model for analyzing user departure time distributions under the UE state, offering insights for potential planning and policy applications, such as fare policies and personalized information provision. Future research will incorporate user heterogeneity into the DTUE-PT model and utilize the model for policy design and evaluation.

\section*{Acknowledgments}
This work was supported by the Monash University for providing the centrally awarded scholarship (Monash Graduate Scholarship and Monash International Tuition Scholarship); the Australian Research Council under Discovery Project [DP190100013]; and the TRENoP (Swedish Strategic Research Area in Transport) Research Center at KTH, Sweden.

\begin{IEEEbiography}[{\includegraphics[width=1in,height=1.25in,clip,keepaspectratio]{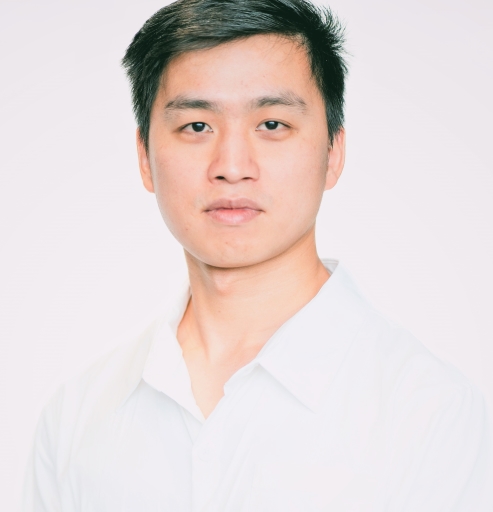}}]{Xia Zhou}
received the B.S. degree in Engineering from Dalian Maritime University, Liaoning, China, and the M.S. degree in engineering from Tianjin University, Tianjin, China. He is pursuing the
Ph.D. degree in the Faculty of Information Technology, Monash University, Melbourne, Australia. His research interests span different techniques and algorithms for optimization, simulation, modeling, and their integration and application to solving transportation assignment problems.\end{IEEEbiography}
\begin{IEEEbiography}[{\includegraphics[width=1in,height=1.25in,clip,keepaspectratio]{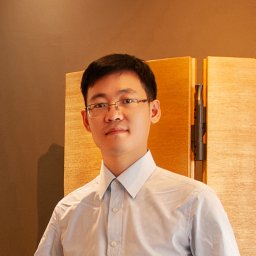}}]{Zhenliang Ma}
received the B.Sc. degree in electrical engineering from Shandong University in 2009, the M.Sc. degree in information technology in 2012, and the Ph.D. degree in transportation engineering from The University of Queensland in 2015. He is currently an Associate Professor of transportation science and a Faculty Member of digital futures with the KTH Royal Institute of Technology. His research interests include statistics, machine learning, computer science-based modeling, simulation, optimization, and control within the framework of selected mobility-related complex systems, which are intelligent transport systems (traffic/public transport/rails) and personal information systems (transport/energy).\end{IEEEbiography}
\begin{IEEEbiography}[{\includegraphics[width=1in,height=1.25in,clip,keepaspectratio]{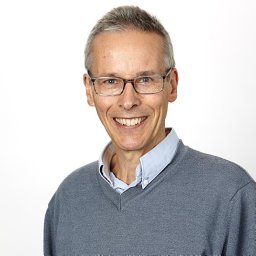}}]{Mark Wallace}
received the bachelor’s degree in mathematics and philosophy from Oxford University, U.K., the master’s degree in artificial intelligence from University of London, U.K., and the
Ph.D. degree from The Southampton University, U.K. He has served as Professor and Associate Dean (Research) with the Faculty ofof Information Technology, Monash University, Melbourne, Australia. His research interests span different techniques and algorithms for optimization and their integration and application to solving complex resource planning and scheduling problems.\end{IEEEbiography}
\begin{IEEEbiography}[{\includegraphics[width=1in,height=1.25in,clip,keepaspectratio]{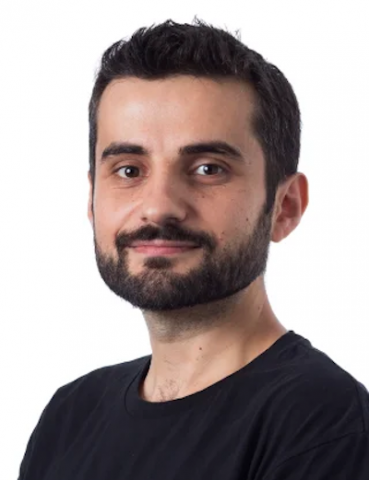}}]{Daniel D. Harabor}
is currently an Associate Professor in the Faculty of Information Technology, Monash University, Melbourne, Australia. He received the B.S. degree (Honours) and the Ph.D. degree in computer science from the Australian National University, Canberra, Australia. His research interests include artificial intelligence, heuristic search, optimisation. In the context of transportation and logistics, he focuses on rail scheduling and capacity analysis of rail supply chains, simulation-based analysis of rail and port operations, data-driven efficacy analysis, and delivery truck fleet optimisation. \end{IEEEbiography}
\begin{appendices}

\section*{Appendix-A: Flow conservation in a congested PT network (Departure time only)} \label{appendix:B}
\renewcommand{\theequation}{A.\arabic{equation}} 
\renewcommand{\thefigure}{A.\arabic{figure}}     
\renewcommand{\thetable}{A.\arabic{table}}       
\setcounter{equation}{0} 
\setcounter{figure}{0}   
\setcounter{table}{0}    

This appendix section describes the model used to evaluate the arrival times of all users in the public transport system. The input and output are $\mathbf{\textit{Q}}$ and  $\mathbf{\textit{A}}$. In this paper, we propose a simulation model to implement Eqs. (\ref{eq:7})-(\ref{eq:29}). All notations are summarized in Table \ref{tab2}.
\subsection*{(1) Picking up users at origin stations}\
$g_o(t)$ consists of the denied boarding users and new arrivals:
\begin{equation}
g_o(t) = h_o(t_l) + \sum_{k \in \{o\} \times D} q_k(t)
\label{eq:7}
\end{equation}

$f_o(t)$ is constrained by train capacity:
\begin{equation}
f_o(t) = \min \left(g_o(t), \Omega_{o(t)}\right)
\label{eq:88}
\end{equation}

To guarantee the flow consistency: 
\begin{equation}
f_o(t) = \sum_{k \in \{o\} \times D}  f_o^k(t)
\label{eq:999}
\end{equation}
\begin{equation}
\sum_{t_1 \in T_o(t)} f_o^k(t_1) \leq \sum_{t_2 \in T_o(t)} q_k(t_2), \quad k \in \{o\} \times {D}
\label{eq:10}
\end{equation}

$h_o(t)$ is updated by the difference between $\Omega_{o(t)}$ and $g_o(t)$:
\begin{equation}
h_o(t) =  \begin{cases}
    0, & g_o(t) <=  \Omega_{o(t)}\\
    g_o(t) - \Omega_{o(t)}, & g_o(t) >  \Omega_{o(t)}
\end{cases}
\label{eq:11}
\end{equation}
\begin{equation}
h_o(t) = \sum_{k \in \{o\} \times D} h_o^k(t)
\label{eq:12}
\end{equation}

To guarantee the flow consistency:
\begin{equation}
h_o^k(t) = h_o^k(t_l) + q_{k}(t) - f_o^k(t), \quad k \in \{o\} \times D 
\label{eq:13}
\end{equation}
\begin{equation}
z_o(t) = f_o(t)
\label{eq:14}
\end{equation}
\subsection*{(2) Dropping off users and picking up users at transfer stations}\

\textit{\textbf{Step 1:}} \textit{Dropping off users}.

$m_b(t)$ consists of various users who need transfer here:
\begin{equation}
m_b(t) = \sum_{k \in K}  m_b^k(t)
\label{eq:15}
\end{equation}
\begin{equation}
m_{b}^k(t) = \sum_{p \in P_{b}(t)} \delta_{b}^k(t) f_{p}^k(t)
\label{eq:16}
\end{equation}

\textbf{\textit{Step 2}}: \textit{Picking up users}.

Similar to $f_{o}(t)$, $f_{b}(t)$ is constrained by available capacity: 
\begin{equation}
f_{b}(t) = \min(g_{b}(t), x_{b}(t))
\label{eq:17}
\end{equation}

To guarantee flow consistency:
\begin{equation}
f_{b}(t) = \sum_{k \in K}  f_b^k(t)
\label{eq:18}
\end{equation}
\begin{equation}
\sum_{t_1 \in T_o(t)} f_b^k(t_1) \leq \sum_{t_2 \in T} \delta_b^{k}(T_o(t), t_2) f_{o(k)}^k(t_2)
\label{eq:19}
\end{equation}

$g_{b}(t)$ is updated by denied boarding users and new arrivals:
\begin{equation}
g_{b}(t) = h_{b}(t_l) + \sum_{t_2 \in T} \sum_{k \in K}  \delta_b^{k} (t, t_2) m_b^k(t_2)
\label{eq:20}
\end{equation}

Available capacity:
\begin{equation}
x_{b}(t) = \Omega_{o(t)} - y_{b}(t)
\label{eq:21}
\end{equation}

$y_{b}(t)$ can be the gap between successfully boarded users and total alighting users:
\begin{equation}
y_{b}(t) =  \sum_{p \in P_{b}(t)} \left(f_p(t) - m_p(t) \right) - m_b(t) 
\label{eq:22}
\end{equation}

Then, $h_b(t)$ can be updated by:
\begin{equation}
h_b(t) =  \begin{cases}
    0, & g_b(t) <=  x_b(t)\\
    g_b(t) - x_b(t), & g_b(t) >  x_b(t)
\end{cases}
\label{eq:23}
\end{equation}
\begin{equation}
h_b(t) =  \sum_{k \in K}h_b^k(t)
\label{eq:24}
\end{equation}
\begin{equation}
h_b^k(t) = h_b^k(t_l) + \sum_{t_2 \in T} \delta_{b}^{k}(t, t_2) m_b^k(t_2) - f_b^k(t) 
\label{eq:399}
\end{equation}
\begin{equation}
z_b(t) = f_b(t) + x_b(t)
\label{eq:25}
\end{equation}
\subsection*{(3) Dropping off users at destination stations}\

To guarantee flow consistency, $m_d(t)$ is given by:
\begin{equation}
m_d(t) = \sum_{k \in O \times \{d\}}  m_d^k(t)
\label{eq:26}
\end{equation}
\begin{equation}
m_d^k(t) = \sum_{p \in P_d(t)} f_p^k(t), \quad k \in O \times \{d\} 
\label{eq:27}
\end{equation}
\begin{equation}
m_d(t) = z_{d_l}(t)
\label{eq:28}
\end{equation}

All users should arrive at their destinations:
\begin{equation}
\sum_{t_1 \in T_{o(k)}}  m_d^k(t_1) = Q_k, \quad k \in O \times \{d\}
\label{eq:29}
\end{equation}
\subsection*{Method for evaluating travel costs of users}\
Steps for the $c_{k}^i (t)$ are as follows: 
\begin{itemize}
    \item {Step 1: Prepare the input $\mathbf{\textit{Q}}$ for simulation};
    \item {Step 2: Evaluate $\mathbf{\textit{M}}$ using the simulation};
    \item {Step 3: Output $\mathbf{\textit{A}}$ based on $\mathbf{\textit{M}}$}: 
    \item {Step 4: Calculate $\mathbf{\textit{C}}$ using Eqs. (\ref{eq:4}) - (\ref{eq:6}) based on $\mathbf{\textit{M}}$}:
\end{itemize}

\section*{Appendix-B: Flow conservation in a congested PT network (Simultaneous Departure time and Route choice)} \label{appendix:B2}
\renewcommand{\theequation}{B.\arabic{equation}} 
\renewcommand{\thefigure}{B.\arabic{figure}}     
\renewcommand{\thetable}{B.\arabic{table}}       
\setcounter{equation}{0} 
\setcounter{figure}{0}   
\setcounter{table}{0}    
In this section, we show the extended model incorporating route choice. Eqs. \ref{eq:7}-\ref{eq:29} in Appendix-A is updated by Eqs. \ref{eq:5555}-\ref{eq:2828} as follows:

Let
\begin{itemize}
    \item $R_k$ denote the set of routes within OD pair $k$;
    \item $(t, r)$ denote the option for departure time $t$ and route $r$;
    \item $(k, t, r, i)$ denote the $ith$ passenger in the group $I_k(t, r)$;
    \item $\delta_b^k(t_1, t_2, r)$: a binary variable that equals 1 when train $t_1$ is the closest arriving train connecting to $t_2$ at $b$ and user group ($k, t_2, r$) need to transfer at $b$; otherwise, 0.
\end{itemize}

\subsubsection{Picking up passengers at origin nodes}\
\begin{equation}
g_o(t) = \sum_{k \in \{o\} \times D} \sum_{r \in R_k} q_k(t, r) + h_o(t_l)
\label{eq:5555}
\end{equation}
\begin{equation}
f_o(t) = \min \left(g_o(t), S_t\right)
\label{eq:6666}
\end{equation} 
\begin{equation}
f_o(t) = \sum_{k \in \{o\} \times D} \sum_{r \in R_k} f_o^k(t, r)
\label{eq:7777}
\end{equation}
\begin{equation}
\sum_{t_1 \in T_o(t)} f_o^k(t_1, r) \leq \sum_{t_2 \in T_o(t)} q_k(t_2, r)
\label{eq:8888}
\end{equation}
\begin{equation}
h_o(t) =  \begin{cases}
    0, & g_o(t) <=  S_t\\
    g_o(t) - S_t, & g_o(t) >  S_t
\end{cases}
\label{eq:9}
\end{equation}
\begin{equation}
h_o(t) = \sum_{k \in \{o\} \times D}\sum_{r \in R_k} h_o^k(t, r)
\label{eq:1010}
\end{equation}
\begin{equation}
h_o^k(t, r) = h_o^k(t_l, r) + q_{k}(t, r) - f_o^k(t, r)
\label{eq:11111}
\end{equation}
\begin{equation}
z_o(t) = f_o(t)
\label{eq:1212}
\end{equation}
\subsubsection{Dropping off passengers and picking up passengers at transfer nodes}\ 

\textit{\textbf{Step 1:}} Dropping off passengers. 
\begin{equation}
m_b(t) = \sum_{k \in K} \sum_{r \in R_k} m_b^k(t, r)
\label{eq:1313}
\end{equation}
\begin{equation}
    m_b^k(t, r) = \delta_b^k {(t, r)} \sum_{p \in P_b(t)} f_p^k(t, r)
\label{eq:1414}
\end{equation}

\textbf{\textit{Step 2}}: Picking up passengers.
\begin{equation}
f_{b}(t) = \min(g_{b}(t), x_{b}(t))
\label{eq:1515}
\end{equation}
\begin{equation}
f_{b}(t) = \sum_{k \in K} \sum_{r \in R_k} f_b^k(t, r)
\label{eq:1616}
\end{equation}
\begin{equation}
\sum_{t_1 \in T_o(t)} f_b^k(t_1, r) \leq \sum_{t_2 \in T} \delta_b^k(T_o(t), t_2, r) f_{o(k)}^k(t_2, r)
\label{eq:1717}
\end{equation}
\begin{equation}
g_{b}(t) = h_{b}(t_l) + \sum_{t_2 \in T} \sum_{k \in K} \sum_{r \in R_k} \delta_b^k(t, t_2, r) m_b^k(t_2, r)
\label{eq:1818}
\end{equation}
\begin{equation}
x_{b}(t) = S_t - y_{b}(t)
\label{eq:1919}
\end{equation}
\begin{equation}
y_{b}(t) =  \sum_{p \in P_{b}(t)} \left(f_p(t) - m_p(t) \right) - m_b(t) 
\label{eq:2020}
\end{equation}
\begin{equation}
h_b(t) =  \begin{cases}
    0, & g_b(t) \leq  x_b(t)\\
    g_b(t) - x_b(t), & g_b(t) >  x_b(t)
\end{cases}
\label{eq:2121}
\end{equation}
\begin{equation}
h_b(t) =  \sum_{k \in K}\sum_{r \in R_k}h_b^k(t, r)
\label{eq:22000}
\end{equation}
\begin{equation}
h_b^k(t, r) = h_b^k(t_l, r) + \sum_{t_2 \in T}\delta_{b}^{k}(t, t_2, r) m_b^k(t_2, r) - f_b^k(t, r)
\label{eq:2323}
\end{equation} 
\begin{equation}
z_b(t) = f_b(t) + x_b(t)
\label{eq:2424}
\end{equation}

\subsubsection{Dropping off passengers at destination nodes}\
\begin{equation}
m_d(t) = \sum_{k \in O \times \{d\}} \sum_{r \in R_k} m_d^k(t, r)
\label{eq:2525}
\end{equation}
\begin{equation}
m_d^k(t, r) = \sum_{p \in P_d(t)} f_p^k(t, r), \quad k \in O \times \{d\}
\label{eq:2626}
\end{equation}
\begin{equation}
m_d(t) = z_{d_l}(t)
\label{eq:2727}
\end{equation}
\begin{equation}
\sum_{t_1 \in T_d} \sum_{r \in R_k} m_d^k(t_1,r) = Q_k, \quad k \in O \times \{d\}
\label{eq:2828}
\end{equation}
\section*{Appendix-C: Proposition 1 proof. } \label{appendix:D}
\renewcommand{\theequation}{C.\arabic{equation}} 
\renewcommand{\thefigure}{C.\arabic{figure}}     
\renewcommand{\thetable}{C.\arabic{table}}       
\setcounter{equation}{0} 
\setcounter{figure}{0}   
\setcounter{table}{0}    
This appendix section uses a simple example to explain why interacting OD flows can cause a non-zero gap in public transport. In this example, we cannot further reduce the option cost gap by shifting users from the non-optimal option to the optimal option, which implies that the gap cannot be reduced to zero. Therefore, the solution of Eqs. (\ref{eq:33}a) may not exist in the public transport system when considering interacting OD flows.

\begin{figure}[ht]
    \centering
    \includegraphics[width=0.5\textwidth]{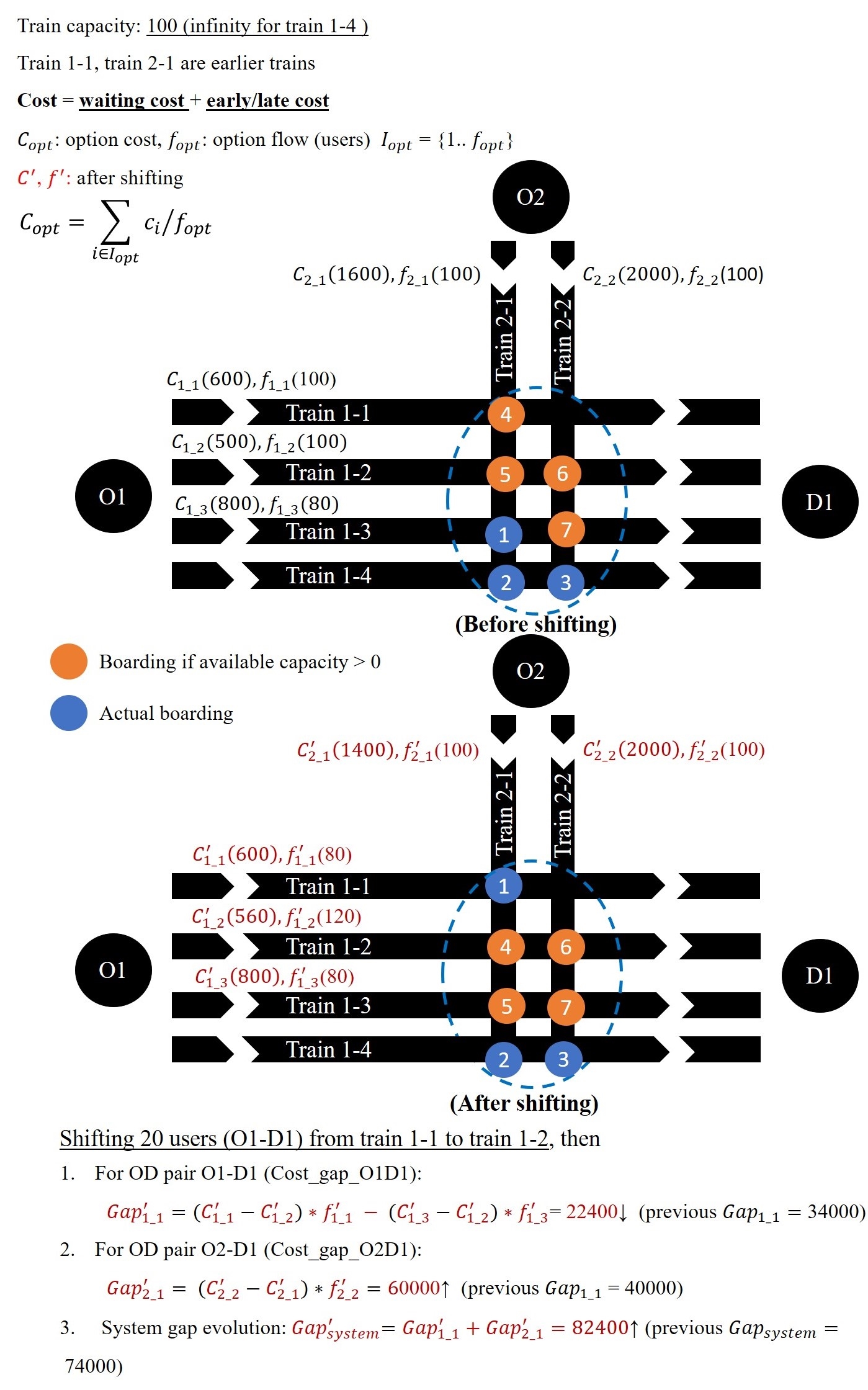}
    \caption{Example of interacting OD flows}
    \label{fig:fig4}
\end{figure}

\textbf{Proof.} Consider the scenario depicted in Fig. \ref{fig:fig4}, where two OD pairs ($O_1D_1$ and $O_2D_1$) share public transport infrastructure. The departure time sequence for trains is as follows: $t_{2\_2}^d > t_{2\_1}^d$, $t_{1\_4}^d > t_{1\_3}^d > t_{1\_2}^d > t_{1\_1}^d$. For users traveling between $O_2$ and $D_1$, the choice lies between train 2-1 or train 2-2, with virtual transfer stations (e.g., station 1, station 2, etc.) facilitating connections between different train lines.

Assuming train 1-2 is the optimal choice for $O_1D_1$ and is at maximum capacity (100), while train 1-1 is full, and train 1-3 has available capacity (20), the initial user distribution over time options is as follows:
\begin{itemize}
\item Option 1-1: 1st to 100th users (all boarding train 1-1)
\item Option 1-2: 101st to 200th users (all boarding train 1-2)
\item Option 1-3: 201st to 280th users (all boarding train 1-3)
\item Option 2-1: 281st to 380th users (281st to 300th users boarding train 1-3, and 301st to 380th users boarding train 1-4 after waiting)
\item Option 2-2: 381st to 480th users (all boarding train 1-4)
\end{itemize}
To reduce the cost gap within OD pair $O_1D_1$, we shift 20 users (81st to 100th) from option 1-1 (choosing train 1-1) to option 1-2 (choosing train 1-2). The new flow distribution is:
\begin{itemize}
\item Option 1-1: 1st to 80th users (all boarding train 1-1)
\item Option 1-2: 81st to 200th users (81st to 100th users boarding train 1-3 after waiting, and 101st to 200th users boarding train 1-2)
\item Option 1-3: 201st to 280th users (all boarding train 1-3)
\item Option 2-1: 281st to 380th users (281st to 300th users boarding train 1-1, and 301st to 380th users boarding train 1-4 after waiting)
\item Option 2-2: 381st to 480th users (all boarding train 1-4)
\end{itemize}

The updated average cost distribution is:
\begin{itemize}
\item Option 1-1: 600 (originally 500)
\item Option 1-2: 560 (waiting cost is 360; 
(500×100+(500+360)×20)/120=560; originally 500)
\item Option 1-3: 800 (originally 800)
\item Option 2-1: 1400 (waiting cost is 1000; (1600×80 + (1600-1000)×20)/100 = 1400; originally 1600)
\item Option 2-2: 2000 (originally 2000)
\end{itemize}

This example shows that reducing the option cost gap for OD pair $O_1D_1$ (from 34,000 to 22,400) results in an increased option cost gap for OD pair $O_2D_1$ (from 40,000 to 60,000), which in turn increases the total system cost gap (from 74,000 to 82,400). This intricate interplay reveals the challenge: in public transport networks with interacting OD flows, achieving a zero travel cost gap for all OD pairs simultaneously is difficult, demonstrating that interacting OD flows can indeed generate a non-zero gap.

\section*{Appendix-D: Sensitivity analysis settings} \label{appendix:E}
\renewcommand{\theequation}{D.\arabic{equation}} 
\renewcommand{\thefigure}{D.\arabic{figure}}     
\renewcommand{\thetable}{D.\arabic{table}}       
\setcounter{equation}{0} 
\setcounter{figure}{0}   
\setcounter{table}{0}    
This appendix section summarizes detailed information about different initial model solution settings for sensitivity analysis. Table \ref{tab:initial-settings} presents these 5 initial model solution settings, where the latest available train is scheduled to depart at 12:25:00 p.m. The initial optimal departure times for OD pairs are detailed in Table \ref{tab:initial-settings}.
\begin{table}[htbp]
  \centering
  \caption{5 Initial model solution settings}
  \label{tab:initial-settings}
  \begin{tabular}{p{0.8cm} p{0.4cm} p{0.52cm} p{0.4cm} p{0.52cm} p{0.4cm} p{0.52cm} p{0.4cm} p{0.52cm}} 
    \toprule
    \textbf{Time} & \multicolumn{4}{c}{\textbf{OD pairs}} & \multicolumn{4}{c}{\textbf{OD pairs (continued)}}\\
    \cmidrule(lr){2-5}
    &  \textbf{1-9} & \textbf{1-10} & \textbf{...} & \textbf{4-12} &  \textbf{1-9} & \textbf{1-10} & \textbf{...} & \textbf{4-12}\\
    & \multicolumn{4}{c}{\textit{Default}} & \multicolumn{4}{c}{\textit{Latest}}\\
    \midrule
    5:00:00 & 0 & 0 & 0 & 0 & 0 & 0 & 0 & 0\\
    5:05:00 & 0 & 0 & 0 & 0 & 0 & 0 & 0 & 0\\
    5:10:00 & 0 & 0 & 0 & 0 & 0 & 0 & 0 & 0\\
    ... & 0 & 0 & 0 & 0 & 0 & 0 & 0 & 0\\
    Preferred & \multicolumn{4}{c}{2000 or 0} & 0 & 0 & 0 & 0\\
    ... & 0 & 0 & 0 & 0 & 0 & 0 & 0 & 0\\
    12:25:00 & 0 & 0 & 0 & 0 & 2000 & 2000 & 2000 & 2000\\
    ... & 0 & 0 & 0 & 0& 0 & 0 & 0 & 0\\
    13:15:00 & 0 & 0 & 0 & 0 & 0 & 0 & 0 & 0\\

    & \multicolumn{4}{c}{\textit{Uniform}} & \multicolumn{4}{c}{\textit{Default-Earliest}}\\
    \midrule
    5:00:00 & 22 & 22 & 22 & 22 & 1000 & 1000 & 1000 & 1000\\
    ... & 22 & 22 & 22 & 22 & 0 & 0 & 0 & 0\\
    Preferred & 22 & 22 & 22 & 22 & \multicolumn{4}{c}{1000 or 0}\\
    ... & 22 & 22 & 22 & 22 & 0 & 0 & 0 & 0\\
    12:25:00 & 22 & 22 & 22 & 22 & 0 & 0 & 0 & 0\\
    12:30:00 & 20 & 20 & 20 & 20 & 0 & 0 & 0 & 0\\
    ... & 0 & 0 & 0 & 0 & 0 & 0 & 0 & 0\\
    13:15:00 & 0 & 0 & 0 & 0 & 0 & 0 & 0 & 0\\
    
    & \multicolumn{4}{c}{\textit{Earliest}} & \multicolumn{4}{c}{\textit{-}}\\
    \midrule
    5:00:00 & 2000 & 2000 & 2000 & 2000 & - & - & - & -\\
    ... & 0 & 0 & 0 & 0 & - & - & - & -\\
    Preferred & 0 & 0 & 0 & 0 & - & - & - & -\\
    ... & 0 & 0 & 0 & 0 & - & - & - & -\\
    12:25:00 & 0 & 0 & 0 & 0 & - & - & - & -\\
    ... & 0 & 0 & 0 & 0 & - & - & - & -\\
    13:15:00 & 0 & 0 & 0 & 0 & - & - & - & -\\
    \bottomrule
  \end{tabular}
\end{table}  

\end{appendices}
\vfill

\end{document}